\documentclass[10pt]{amsart}

\usepackage{amssymb,amsmath,amsthm,color}
\usepackage[all]{xy}

\usepackage[OT2,OT1]{fontenc}
\newcommand\cyr{%
\renewcommand\rmdefault{wncyr}%
\renewcommand\sfdefault{wncyss}%
\renewcommand\encodingdefault{OT2}%
\normalfont \selectfont} \DeclareTextFontCommand{\textcyr}{\cyr}

\newtheorem{theorem}{Theorem}[section]
\newtheorem{proposition}[theorem]{Proposition}
\newtheorem{lemma}[theorem]{Lemma}
\newtheorem{corollary}[theorem]{Corollary}
\theoremstyle{definition} 
\newtheorem{definition}[theorem]{Definition}
\newtheorem{conjecture}[theorem]{Conjecture}
\newtheorem{notation}[theorem]{Notation}

\newtheorem{remark}[theorem]{Remark}

\DeclareMathOperator{\Gal}{Gal}
\DeclareMathOperator{\Real}{Re}
\DeclareMathOperator{\alg}{alg}\DeclareMathOperator{\SL}{SL}
\DeclareMathOperator{\Sp}{Sp}\DeclareMathOperator{\M}{M}
\DeclareMathOperator{\GL}{GL}
\DeclareMathOperator{\modu}{mod}\DeclareMathOperator{\cusp}{cusp}
\DeclareMathOperator{\GCD}{gcd}\DeclareMathOperator{\spin}{spin}
\DeclareMathOperator{\SGN}{sgn}\DeclareMathOperator{\Tr}{Tr}
\DeclareMathOperator{\diag}{diag}\DeclareMathOperator{\st}{st}
\DeclareMathOperator{\cohom}{H}
\DeclareMathOperator{\new}{new}\DeclareMathOperator{\Frob}{Frob}\DeclareMathOperator{\Ann}{Ann}
\DeclareMathOperator{\ev}{ev}\DeclareMathOperator{\cont}{cont}\DeclareMathOperator{\crys}{crys}
\DeclareMathOperator{\Fil}{Fil}\DeclareMathOperator{\ur}{ur}
\DeclareMathOperator{\im}{im}\DeclareMathOperator{\Ext}{Ext}\DeclareMathOperator{\Sha}{\textcyr{Sh}}
\DeclareMathOperator{\vol}{vol}\DeclareMathOperator{\semi}{\, ss}\DeclareMathOperator{\Nm}{Nm}\DeclareMathOperator{\disc}{disc}
\DeclareMathOperator{\image}{Image}\DeclareMathOperator{\J}{J}

\newcommand{\inte}{\mathbb{Z}}
\newcommand{\nat}{\mathbb{N}}
\newcommand{\rat}{\mathbb{Q}}
\newcommand{\real}{\mathbb{R}}
\newcommand{\cmplx}{\mathbb{C}}
\newcommand{\F}{\mathbb{F}}
\newcommand{\T}{\mathbb{T}}
\newcommand{\A}{\mathbb{A}}
\newcommand{\B}{\mathbb{B}}
\newcommand{\q}{\mathfrak{q}}
\newcommand{\Ff}{F_{f}}
\newcommand{\Gone}{\Gamma_{1}}
\newcommand{\Gtwo}{\Gamma_{2}}
\newcommand{\ROI}{\mathcal{O}}
\newcommand{\h}[1]{\mathfrak{h}^{#1}}
\newcommand{\f}{\textbf{f}}
\newcommand{\Slash}[1]{|_{#1}}
\newcommand{\mat}[4]{\begin{pmatrix}{#1}&{#2}\\{#3}&{#4}\end{pmatrix}}
\newcommand{\p}{\mathfrak{p}}
\newcommand{\absgal}{\Gal(\overline{\rat}/\rat)}

\begin{document}

\title{Saito-Kurokawa Lifts and Applications to the Bloch-Kato Conjecture}
\author{Jim Brown}
\address{Department of Mathematics\\
The Ohio State University\\ Columbus, OH  43210}
\email{jimlb@math.ohio-state.edu}
\subjclass[2000]{Primary 11F33, 11F67; Secondary 11F46, 11F80}

\begin{abstract}  
Let $f$ be a newform of weight $2k-2$ and level $1$.  In this paper
we provide evidence for the Bloch-Kato conjecture for modular forms.  We 
demonstrate an implication that under suitable hypotheses if
$\varpi \mid L_{\alg}(k,f)$ then $p \mid \# \cohom_{f}(\rat,W_{f}(1-k))$ where
$p$ is a suitably chosen prime and $\varpi$ a uniformizer of a finite extension
$K/\rat_{p}$.  We demonstrate this by establishing a congruence between the Saito-Kurokawa
lift $\Ff$ of $f$ and a cuspidal Siegel eigenform $G$ that is not a Saito-Kurokawa lift.  We then
examine what this congruence says in terms of Galois representations to produce a non-trivial $p$-torsion
element in $\cohom_{f}^1(\rat,W_{f}(1-k))$.
\end{abstract}

\maketitle

\section{Introduction}
Let $f$ be a newform of weight $2k-2$ and level $1$. The Bloch-Kato conjecture for 
modular forms roughly states that the special values of the $L$-function associated to $f$
should measure the size of the corresponding Selmer groups.  In this paper 
we will demonstrate under suitable hypotheses that if $\varpi \mid L_{\alg}(k,f)$, then $p \mid \#\cohom_{f}^1(\rat,W_{f}(1-k))$
where $p$ is a suitably chosen prime and $\varpi$ is a uniformizer of a finite extension $K/\rat_{p}$.  
For a precise statement see Theorem \ref{thm:finaltheorem}.  

The general outline of the method of proof of Theorem \ref{thm:finaltheorem} goes back to Ribet's proof of the
converse of Herbrand's theorem (\cite{ribet}), which was then extended by Wiles in his proof of the main conjecture of Iwasawa
theory for totally real fields (\cite{wiles}).  The method used by Ribet and Wiles is as follows.  Given a positive integer $k$ and a primitive Dirichlet character $\chi$ 
of conductor $N$ so that $\chi(-1) = (-1)^{k}$, one has an associated Eisenstein series $E_{k,\chi}$.  For a prime $p \nmid N$, one can show that there is a cuspidal eigenform $g$ of weight $k$ and level $M$ with $N \mid M$ so that 
$g \equiv E_{k,\chi} (\modu \p)$ for some prime $\p \mid p$.  This congruence is used to 
study the residual Galois representation of $g$.  It is shown that $\overline{\rho}_{g,\p} \simeq \mat{1}{*}{0}{\chi \omega^{k-1}}$ is non-split
where $\omega$ is the reduction of the $p$-adic cyclotomic character.  This allows one to show that $*$ gives a non-zero cohomology class
in $\cohom_{\ur}^1(\rat, \chi^{-1}\omega^{1-k})$.   

For our purposes, the character in the Ribet/Wiles' method will be replaced with a newform $f$ of weight $2k-2$ and level 1.  
Associated to $f$ we have its
Saito-Kurokawa lift $\Ff$, our replacement for the Eisenstein series $E_{k,\chi}$.  Our goal is to find a cuspidal Siegel
eigenform $G$ that is not a Saito-Kurokawa lift so that the eigenvalues of $G$ are congruent modulo $\varpi$ to those of $\Ff$.  We are
able to produce such a $G$ by exploiting the explicit nature of the Saito-Kurokawa correspondence.  Also central to producing $G$ is an inner
product relation due to Shimura (\cite{shimurazfctns}).  In order to assure the $G$ we construct is not a Saito-Kurokawa lift, we are
forced to act on $G$ with a particular Hecke operator that kills all Saito-Kurokawa lifts other then $\Ff$. It is in this step that
we must insert the hypothesis that $f$ is ordinary at $p$.  It appears this is merely a technical restriction that we hope to remove in a 
subsequent paper.  For the precise statement of the congruence see Theorem \ref{thm:thecongruence}.

Once we have produced a congruence modulo $\varpi$ between the Hecke eigenvalues of $\Ff$ and $G$, we study the
associated 4-dimensional Galois representations.  Again we use the explicit nature of the Saito-Kurokawa correspondence to conclude that
$\overline{\rho}_{\Ff} \simeq \omega^{k-2} \oplus \overline{\rho}_{f} \oplus \omega^{k-1}$.  Using our congruence we are able to determine that 
$\overline{\rho}^{\semi}_{G} \simeq \overline{\rho}_{\Ff}$.  From this we deduce the form of $\overline{\rho}_{G}$ by adapting arguments 
in \cite{ribet} to the 4-dimensional case and applying results of \cite{skinnerurban} on the necessary shape of $\rho_{G}$.  Some elementary 
arguments using class field theory allow us to conclude that we have a non-zero torsion element of the Selmer group $\cohom_{f}^1(\rat,W_{f}(1-k))$.
We conclude with a non-trivial numerical example of Theorem \ref{thm:finaltheorem} with $p = 516223$ and $f$ of weight $54$. 

While this paper only deals with the case of full level, it is anticipated that similar results hold true for odd square-free level. We hope to treat 
the case of odd square-free level in a subsequent paper.

The author would like to thank Chris Skinner for many helpful conversations.

\section{Notation and definitions}\label{sec:notation}

In this section we fix notation and definitions that will be used throughout this 
paper.  

Denote the adeles over $\rat$ by $\A$.  We let $\textbf{f}$ denote the finite set of
places.  For $p$ a prime number, we fix once and for all compatible embeddings $\overline{\rat} \hookrightarrow \overline{\rat}_{p}$,
$\overline{\rat} \hookrightarrow \cmplx$, and $\overline{\rat}_{p} \hookrightarrow \cmplx$. Let $\varepsilon_{p}$ be the $p$-adic
cyclotomic character $\varepsilon_{p}: \absgal \rightarrow \GL_1(\inte_{p})$.  Recall that $\varepsilon_{p}$ is unramified away from $p$
and one has $\varepsilon_{p}(\Frob_{\ell}) = \ell$ for $\ell \neq p$.  We write $\rat_{p}(n)$ for the 1-dimensional space over $\rat_{p}$
on which $\absgal$ acts by $\varepsilon_{p}^{n}$ and similarly for $\inte_{p}(n)$.  We denote the residual representation of 
$\varepsilon_{p}$ by $\omega_{p}$.  We will drop the $p$ when it is clear from the context.

Let $\Sigma$ be a set of primes.  For an $L$-function we write $L^{\Sigma}$ to denote the restricted Euler product of $L$ 
over primes not in $\Sigma$ and $L_{\Sigma}$ to denote the restricted Euler product over primes in $\Sigma$.

For a ring $R$, we let $\M_{n}(R)$ denote the set of $n$ by $n$ matrices with entries in $R$. For a matrix $x \in \M_{2n}(R)$, we write
\begin{equation*}
x = \mat{a_{x}}{b_{x}}{c_{x}}{d_{x}}
\end{equation*}
where $a_{x}$, $b_{x}$, $c_{x}$, and $d_{x}$ are all in $\M_{n}(R)$.  We drop the subscript $x$ when it is clear from the context.

Denote the group $\SL_2(\inte)$ by $\Gamma_1$.  We refer to a subgroup of $\Gamma_1$ as a congruence subgroup if it contains $\Gamma(N)$ for 
some positive integer $N$.  We denote the complex upper half-plane by $\h{1}$.  As usual, $\GL_2^{+}(\real)$ acts on 
$\h{1} \cup \real \cup \{ \infty \}$ via linear fractional transformations.  We let $\Gamma_1^{\J} = \Gamma_1 \ltimes \inte^{2}$ be the full Jacobi
modular group, as defined in \cite{eichlerzagier}.  Recall that the symplectic group is defined by 
\begin{equation*}
\Sp_{2n}(\real) = \{ \gamma \in \M_{2n}(\real) :\, ^{t}\!\gamma\iota_{n} \gamma = \iota_{n}\}, \quad \iota_{n} = \mat{0_{n}}{-1_{n}}{1_{n}}{0_{n}}
\end{equation*}
where we write $1_{n}$ to denote the $n$ by $n$ identity matrix.
We denote $\Sp_{2n}(\inte)$ by $\Gamma_{n}$.  Siegel upper half-space is given by 
\begin{equation*}
\h{n} = \{ Z \in \M_{n}(\cmplx):\, ^{t}\!Z = Z,\, \text{Im}(Z) >0 \}.
\end{equation*}
Siegel upper half-space comes equipped with an action of $\Sp_{2n}(\real)$ given by 
\begin{equation*}
\mat{A}{B}{C}{D} Z = (AZ + B) (CZ+D)^{-1}.
\end{equation*}

For a congruence subgroup $\Gamma \subseteq \Gamma_1$, we write $M_{k}(\Gamma)$ to denote the
space of modular forms of weight $k$ on the congruence subgroup $\Gamma$.  For $f \in M_{k}(\Gamma)$, we denote
the $n^{\text{th}}$ Fourier coefficient of $f$ by $a_{f}(n)$.  Given a ring $R \subseteq \cmplx$, we write $M_{k}(\Gamma,R)$ 
to denote the space of modular forms with Fourier coefficients in $R$.  Let $S_{k}(\Gamma)$ denote
the space of cusp forms.  For $f_1,f_2 \in M_{k}(\Gamma)$ with $f_1$ or $f_2$ a cusp form, the Petersson product is given by
\begin{equation*}
\langle f_1, f_2 \rangle = \frac{1}{[\overline{\Gamma}_1:\overline{\Gamma}]} \int_{\Gamma\backslash \h{1}} f_1(z) \overline{f_2(z)} y^{k-2}dx\, dy
\end{equation*}
where $\overline{\Gamma}_1$ means $\Gamma_1/\pm 1_2$ and $\overline{\Gamma}$ is the image of $\Gamma$ in $\overline{\Gamma}_1$. 
We write $\T_{R}(\Gamma)$ for the usual Hecke algebra over the ring $R$ for the congruence subgroup $\Gamma$.  We drop $\Gamma$ 
from the notation when it is clear from the context.  We say $f$ is a newform 
if it is an eigenform for all the Hecke operators $T(n)$ with Fourier expansion normalized so that the Fourier coefficients are equal
to the eigenvalues. We write $S_{k}^{\new}(\Gamma)$ to denote the space of newforms.  

The only half-integral weight modular forms we will be interested in are the ones in Kohnen's $+$-space defined by
\begin{equation*}
S_{k-1/2}^{+}(\Gamma_0(4)) = \{ g \in S_{k-1/2}(\Gamma_0(4)):\, \text{$a_{g}(n) = 0$ if $(-1)^{k-1}n \equiv 2,3 (\modu 4)$}\}.
\end{equation*}
The Petersson product on $S_{k-1/2}^{+}(\Gamma_0(4))$ is given by
\begin{equation*}
\langle g_1, g_2 \rangle = \int_{\Gamma_0(4)\backslash \h{1}} g_1(z) \overline{g_2(z)} y^{k-5/2} dx\, dy.
\end{equation*}

We denote the space of Jacobi cusp forms on $\Gamma_1^{\J}$ by $J_{k,1}^{\cusp}(\Gamma_1^{\J})$.  The inner product is given
by 
\begin{equation*}
\langle \phi_1, \phi_2 \rangle = \int_{\Gamma_1^{\J}\backslash \h{1} \times \cmplx} \phi_1(\tau, z)\overline{\phi_2(\tau,z)} v^{k-3} e^{-4 \pi y^2/v} dx \,dy\, du\, dv
\end{equation*}
for $\phi_1, \phi_2 \in J_{k,1}^{\cusp}(\Gamma_1^{\J})$ and $\tau= u + i v$, $z = x+iy$. 

Given a congruence group $\Gamma \subseteq \Sp_{2n}(\inte)$, we denote the space of Siegel modular forms of weight $k$ for $\Gamma$ by
$\mathcal{M}_{k}(\Gamma)$.  The space of cusp forms is denoted by $\mathcal{S}_{k}(\Gamma)$.  For $\gamma \in \Sp_{2n}^{+}(\real)$, the slash operator
of $\gamma$ on a Siegel modular form $F$ of weight $k$ is given by $(F|_{k}\gamma)(Z) = \det(C_{\gamma}Z + D_{\gamma})^{-k}F(\gamma Z)$. For $F$ and $G$ two Siegel modular forms
with at least one of them a cusp form for $\Gamma$ of weight $k$, define the Petersson product of $F$ and $G$ by
\begin{equation*}
\langle F, G \rangle = \frac{1}{[\overline{\Gamma}_{n}:\overline{\Gamma}]} \int_{\Gamma\backslash \h{n}} F(Z) \overline{G(Z)} \det(Y)^{k} d\mu(Z).
\end{equation*}
We write $\T_{S,R}(\Gamma)$ for the usual Hecke algebra generated over $R$ by the Hecke operators on Siegel modular forms for the congruence group $\Gamma$.  We drop $\Gamma$ from the notation when it is clear from the context.  For a thorough treatment of Hecke operators on Siegel modular 
forms one can consult \cite{andrianov}.

We will mainly be interested in the case when $F \in \mathcal{S}_{k}(\Gtwo)$. Let $F \in \mathcal{S}_{k}(\Gtwo)$ be a Hecke eigenform with eigenvalues $\lambda_{F}(m)$.  The standard zeta function
associated to $F$ is given by 
\begin{equation}\label{eqn:standardzeta}
L_{\st}(s,F) = \prod_{\ell} W_{\ell}(\ell^{-s})^{-1}
\end{equation}
where 
\begin{equation*}
W_{\ell}(t) = (1-\ell^{2}t) \prod_{i=1}^{2} (1- \ell^2 \alpha_{\ell,i}t)(1-\ell^2\alpha_{\ell,i}^{-1} t)
\end{equation*}
with $\alpha_{\ell,i}$ denoting the Satake parameters. Given a Hecke character $\phi$, the twisted standard zeta function is given by
\begin{equation*}
L_{\st}(s,F,\phi) = \prod_{\ell} W_{\ell}(\phi(\ell)\ell^{-s})^{-1}.
\end{equation*}
Associated to $F$ is another $L$-function called the Spinor $L$-function.  It is defined by 
\begin{equation*}
L_{\spin}(s,F) = \zeta(2s-2k+4) \sum_{m=1}^{\infty} \lambda_{F}(m)m^{-s}.
\end{equation*}
We will also be interested in the Maass space $\mathcal{M}_{k}^{*}(\Gtwo) \subset \mathcal{M}_{k}(\Gtwo)$.  A Siegel modular form $F$ is in the Maass space if the Fourier coefficients of $F$
satisfy the relation
\begin{equation*}
A_{F}(n,r,m) = \sum_{d \mid \GCD(n,r,m)} d^{k-1} A_{F}\left(\frac{nm}{d^2}, \frac{r}{d},1\right)
\end{equation*}
for every $m,n,r \in \inte$ with $m,n,4mn-r^2 \geq 0$ (\cite{zagier}).

\section{The Saito-Kurokawa correspondence}\label{sec:sk}

In this section we review the explicit formula approach to the Saito-Kurokawa correspondence  
established by Maass (\cite{maass1} - \cite{maass3}), Andrianov \cite{andrianovsk}, 
and Zagier \cite{zagier}.  We do not claim a complete account and are mainly concerned with stating 
the relevant facts we need in this paper.  The interested reader is urged to consult the references for 
the details. 

\subsection{The correspondence}
The first step in establishing the Saito-Kurokawa correspondence is to relate the integer weight cusp forms
of weight $2k-2$ and level 1 to half-integer weight modular forms of weight $k-1/2$ and level 4.   This is accomplished
via the Shimura and Shintani liftings.  These maps are adjoint on cusp forms with respect
to the Petersson products. Let $D$ be a fundamental discriminant with $(-1)^{k-1}D >0$.  The Shimura 
lifting $\zeta_{D}$ is a map from $S_{k-1/2}^{+}(\Gamma_0(4))$ to $S_{2k-2}(\Gone)$.  Explicitly, for 
\begin{equation*}
g(z) = \sum c_{g}(n) q^{n} \in S_{k-1/2}^{+}(\Gamma_0(4M))
\end{equation*}
one has
\begin{equation*}
\zeta_{D}g(z) = \sum_{n = 1}^{\infty}\left( \sum_{d \mid n} \left( \frac{D}{d}\right) d^{k-2} c_{g}(|D|n^2/d^2)\right)q^{n}
\end{equation*}
where the summation defining $g(z)$ is over all $n \geq 1$ so that $(-1)^{k-1}n \equiv 0,1(\modu 4)$. 
On the other hand, the Shintani lifting $\zeta_{D}^{*}$ is a map from $S_{2k-2}(\Gone)$ to $S_{k-1/2}^{+}(\Gamma_0(4))$.  One can 
consult \cite{kohnenshimura} for a precise defintion of the Shintani map as its precise definition will not be needed here.  Using these liftings, one
has the following theorem:

\begin{theorem} (\cite{kohnenmathann}) For $D$ a fundamental discriminant with $(-1)^{k-1} D >0$, the Shimura and Shintani liftings give  Hecke-equivariant 
isomorphisms between $S_{2k-2}(\Gone)$ and $S_{k-1/2}^{+}(\Gamma_0(4))$.
\end{theorem}
 
Let $\ROI$ be a ring so that an embedding of $\ROI$ into $\cmplx$ exists.  Choose such an embedding and identify
$\ROI$ with its image in $\cmplx$ via this embedding.  Assume that $\ROI$ contains all the Fourier coefficients of $f$.  The 
Shintani lifting $g_{f}:=\zeta_{D}^{*}f$ is determined only up to normalization by a constant multiple.  However, we do have the following result of Stevens.

\begin{theorem}\label{thm:stevenstheorem} (\cite{stevens}, Prop. 2.3.1) Let $f \in S_{2k-2}(\Gone)$ be a newform. If the Fourier 
coefficients of $f$ are in $\ROI$ then there exists a corresponding Shintani lifting $g_{f}$ of $f$ with Fourier 
coefficients in $\ROI$ as well.
\end{theorem}

\begin{remark} Throughout this paper we fix our $g_{f}$ to have Fourier coefficients in $\ROI$ as in 
Theorem \ref{thm:stevenstheorem}.  If, in addition, $\ROI$ is a discrete valuation ring, we fix our $g_{f}$ to have Fourier 
coefficients in $\ROI$ with some Fourier coefficient in $\ROI^{\times}$.
\end{remark}

We have the following theorem relating half-integral weight cusp forms to Jacobi forms.

\begin{theorem}(\cite{eichlerzagier}, Theorem 5.4) The map defined by 
\begin{equation*}
\sum_{\small\begin{array}{c} D<0, r \in \inte \\ D \equiv r^2 (\modu 4)\end{array}} \hspace*{-.1in} c(D,r) e\left(\frac{r^2-D}{4}\tau + rz \right) \mapsto
\hspace*{-.2in}\sum_{\small \begin{array}{c} D<0\\ D \equiv 0, 1(\modu 4)\end{array}} \hspace*{-.2in} c(D) e(|D|\tau)
\end{equation*}
is a canonical Hecke equivariant isomorphism between $J_{k,1}^{\cusp}(\Gone^{\J})$ and $S_{k-1/2}^{+}(\Gamma_0(4))$ preserving the Hilbert
space structures.
\end{theorem}

Our final step is to relate Jacobi forms to Siegel forms.  Let $F \in \mathcal{S}_{k}^{*}(\Gtwo)$.  One has that $F$ admits a 
Fourier-Jacobi expansion
\begin{equation*}
F(\tau,z,\tau') = \sum_{m \geq 0} \phi_{m}(\tau,z)e(m\tau')
\end{equation*}
where the $\phi_{m}$ are Jacobi forms of weight $k$, index $m$, and level 1.  

\begin{theorem}\label{thm:jacobitosiegel} (\cite{eichlerzagier}, Theorem 6.2) The association $F \mapsto \phi_1$ gives a Hecke equivariant 
isomorphism between $\mathcal{S}_{k}^{*}(\Gtwo)$ and $J_{k,1}^{\cusp}(\Gone^{\J})$.  The inverse map is given by sending 
$\phi(\tau,z) \in J_{k,1}^{\cusp}(\Gone^{\J})$ to $\displaystyle F(\tau,z,\tau') = \sum_{m \geq 0} V_{m} \phi(\tau,z) e(m\tau')$
where $V_{m}$ is the index shifting operator as defined in (\cite{eichlerzagier}, Section 4).
\end{theorem}

\begin{corollary} Let $\phi \in J_{k,1}^{\cusp}(\Gone^{\J},\ROI)$ where $\ROI$ is some ring.  If $F$ is the Siegel modular form 
associated to $\phi$ in Theorem \ref{thm:jacobitosiegel} then $F$ has Fourier coefficients in $\ROI$.
\end{corollary}

\begin{proof} Using that $F$ is in the Maass space and the definition of $V_{m}$ we obtain
\begin{equation*}
A(n,r,m) = \hspace*{-.1in}\sum_{d \mid \GCD(m,n,r)}\hspace*{-.2in} d^{k-1} c\left( \frac{4nm-r^2}{d^2}, \frac{r}{d}\right)
\end{equation*}
where the $c(D,r)$ are the Fourier coefficients of $\phi$.  The rest is clear.
\end{proof}

\noindent Combining these results one obtains the Saito-Kurokawa correspondence.

\begin{theorem}\label{thm:skcorrespondence} (\cite{zagier}) The space $\mathcal{S}_{k}^{*}(\Gtwo)$ is spanned
by Hecke eigenforms.  These are in 1-1 correspondence with newforms $f \in S_{2k-2}(\Gone)$, the correspondence
being such that if $\Ff$ correponds to $f$, then one has
\begin{equation}\label{eqn:skequation}
L_{\spin}(s,\Ff) = \zeta(s-k+1)\zeta(s-k+2)L(s,f).
\end{equation}
\end{theorem}

\begin{corollary}\label{thm:skfourier} The Saito-Kurokawa isomorphism is a Hecke-equivariant isomorphism over $\ROI$. In particular, if $\ROI$ is a discrete valuation ring, 
$\Ff$ has a Fourier coefficient in $\ROI^{\times}$.
\end{corollary}

We also note the following theorem giving an equation relating $\langle \Ff, \Ff \rangle$ to $\langle f, f \rangle$. 

\begin{theorem}\label{thm:innerproducts} (\cite{kohnenskoruppa}, \cite{kohnenzagier}) Let $f \in S_{2k-2}(\Gone)$ be a newform, $\Ff \in \mathcal{S}_{k}^{*}(\Gtwo)$
the corresponding Saito-Kurokawa lift, and $g(z) = \sum c_{g}(n)q^{n}$ the weight $k-1/2$ cusp form corresponding
to $f$ under the Shintani map.  We have the following inner product relation
\begin{equation*}
\langle \Ff, \Ff \rangle = \frac{(k-1)}{2^5 3^2 \pi }\, \cdot \frac{ c_{g}(|D|)^2}{ |D|^{k-3/2}}\cdot \, \frac{L(k,f)}{L(k-1,f,\chi_{D})} \, \langle f, f \rangle
\end{equation*}
where $D$ is a fundamental discriminant so that $(-1)^{k-1} D > 0$ and $\chi_{D}$ is the quadratic character associated to $D$.
\end{theorem}

The standard zeta function of $\Ff$ can be factored into a particularly simple form, as given in the following theorem.

\begin{theorem}\label{thm:standardzetafactorization} Let $N$ be a positive integer, $\Sigma$ the set of primes dividing $N$, and $\chi$ a 
Dirichlet character of conductor $N$.  
Let $f \in S_{2k-2}(\Gone)$ be a newform and $\Ff$ the corresponding Saito-Kurokawa lift of $f$.  The standard zeta function of $\Ff$ factors 
as 
\begin{equation*}
L^{\Sigma}_{\st}(2s, \Ff, \chi) = L^{\Sigma}(2s-2,\chi) L^{\Sigma}(2s+k-3,f,\chi)L^{\Sigma}(2s+k-4,f,\chi).
\end{equation*}
\end{theorem}

\begin{proof} To prove this theorem we need to relate the Satake parameters $\alpha_{i}:=\alpha_{p,i}$ to the
eigenvalues of $f$ in order to decompose the standard zeta function. To accomplish this, we use the following formula
(see \cite{panchishkin}):
\begin{equation*}
L_{\spin, (p)}(s,\Ff) = (1-\alpha_{0}p^{-s})
(1-\alpha_{0}\alpha_1 p^{-s})(1-\alpha_0 \alpha_2 p^{-s})
(1-\alpha_0 \alpha_1 \alpha_2 p^{-s}).
\end{equation*}
Recall that by Equation \ref{eqn:skequation} we have
\begin{equation*}
L_{\spin, (p)}(s,\Ff) =
(1-p^{k-1-s})(1-p^{k-s-2})(1-a_{f}(p)p^{-s} + p^{2k-3-2s}).
\end{equation*}
Letting $x=p^{-s}$, we have the polynomial identity
\begin{equation*}
(1-\alpha_{0}x) (1-\alpha_{0}\alpha_1 x)(1-\alpha_0 \alpha_2 x)
(1-\alpha_0 \alpha_1 \alpha_2
x)=(1-p^{k-1}x)(1-p^{k-2}x)(1-a_{f}(p)x + p^{2k-3}x^2).
\end{equation*}
Therefore we have that $\left \{\alpha_0, \alpha_0 \alpha_1, \alpha_0
\alpha_2, \alpha_0 \alpha_1 \alpha_2 \right \} = \left\{ p^{k-1},
p^{k-2}, \frac{2p^{2k-3}}{a_{f}(p) \pm \sqrt{a_{f}(p)^2 -
4p^{2k-3}}} \right\}$. The values $\alpha_0 \alpha_1$ and
$\alpha_0 \alpha_2$ are completely symmetrical so we set
\begin{equation*}
\alpha_0 \alpha_1 = \frac{2p^{2k-3}}{a_{f}(p) + \sqrt{a_{f}(p)^2 -
4p^{2k-3}}}
\end{equation*}
and
\begin{equation*}
\alpha_0 \alpha_2 = \frac{2p^{2k-3}}{a_{f}(p) - \sqrt{a_{f}(p)^2 -
4p^{2k-3}}}.
\end{equation*}
Since we have $\alpha_0^2 \alpha_1 \alpha_2 = p^{2k-3}$, $\alpha_0
=p^{k-1}$ or $p^{k-2}$ but is arbitrary up to this choice. We fix
$\alpha_0 = p^{k-1}$. Pick $\alpha_{p}$ and $\beta_{p}$ such that
\begin{equation*}
\alpha_{p} + \beta_{p} = a_{f}(p)
\end{equation*}
and
\begin{equation*}
\alpha_{p} \beta_{p} = p^{2k-3}.
\end{equation*}
Thus,
\begin{equation*}
\alpha_1 = \beta_{p} p^{1-k}
\end{equation*}
and
\begin{equation*}
\alpha_2 = \alpha_{p} p^{1-k}.
\end{equation*}
Therefore we can write
\begin{equation*}
(1-\chi(p)\alpha_1 p^{2-2s})(1-\chi(p)\alpha_2 p^{2-2s})=
1-\chi(p)a_{f}(p)p^{3-2s-k} + \chi(p)^2 p^{3-4s}
\end{equation*}
and
\begin{equation*}
(1-\chi(p)\alpha_1^{-1} p^{2-2s})(1-\chi(p)\alpha_2^{-1}
p^{2-2s})= 1-\chi(p)a_{f}(p) p^{4-2s-k} + \chi(p)^2 p^{5-4s}.
\end{equation*}
Substituting this back in for $L^{\Sigma}(2s,\Ff,\chi)$ we
have the result.
\end{proof}

\section{Eisenstein series}\label{sec:eseries}

In this section we study an Eisenstein series $E(Z,s,\chi)$ as
defined by Shimura (\cite{shimuraeseries}, \cite{shimuraeprods},
\cite{shimurabook}).  We begin with basic definitions and then
move to a study of the Fourier coefficients of the Eisenstein
series. We show that under a suitable normalization for a
certain value of $s$ we have that $E(Z,s,\chi)$ is a holomorphic
Siegel modular form with Fourier coefficients that are
$p$-integral for a prime $p$ of our choosing. We next move to
studying an inner product relation of Shimura that calculates the
inner product of $E(Z,s,\chi)$ with a Siegel cusp form $F$ in
terms of $F$ and the standard zeta function associated to $F$.

\subsection{Basic definitions}

Before we can define the Eisenstein series we need to define some
subgroups of $\Sp_{2n}(\A)$ and $\Sp_{2n}(\rat)$.  Let
$\mathfrak{a}$ and $\mathfrak{b}$ be non-zero ideals in $\inte$.
Set
\begin{equation*}
D[\mathfrak{a},\mathfrak{b}] = \Sp_{2n}(\real) \prod_{\ell \in
\textbf{f}} D_{\ell}[\mathfrak{a},\mathfrak{b}]
\end{equation*}
where
\begin{equation*}
D_{\ell}[\mathfrak{a},\mathfrak{b}] = \left\{ x \in
\Sp_{2n}(\rat_{\ell}) : a_{x} \in \M_{n}(\inte_{\ell}), b_{x} \in
\M_{n}(\mathfrak{a}_{\ell}), c_{x} \in
\M_{n}(\mathfrak{b}_{\ell}), d_{x} \in \M_{n}(\inte_{\ell}) \right\}.
\end{equation*}
Define a maximal compact subgroup $C_{\upsilon}$ of
$\Sp_{2n}(\rat_{\upsilon})$ by
\begin{equation*}
C_{\upsilon} = \left\{ \begin{array}{ll} \{\alpha \in
\Sp_{2n}(\real) : \alpha (i) = i\} & \upsilon = \infty,\\
\Sp_{2n}(\rat_{\upsilon}) \cap \GL_{2n}(\inte_{\upsilon}) &
\upsilon \in \textbf{f}, \end{array} \right.
\end{equation*}
and set $\displaystyle C = \prod C_{\upsilon}$.  It is understood here that $i$ denotes the $n \times n$
identity matrix multiplied by the complex number $i$. Let $P$ be the Siegel parabolic of
$\Sp_{2n}(\rat)$ defined by
\begin{equation*}
P = \left\{ x \in \Sp_{2n}(\rat) : c_{x} = 0 \right\}.
\end{equation*}
Set
\begin{equation*}
\mathbb{S}^{n}(R) = \{ x \in \M_{n}(R) : \, ^{t}x = x \}.
\end{equation*}
We write elements $Z \in \h{n}$ as $Z = X + i Y$
with $X, Y \in \mathbb{S}^{n}(\real)$ and $Y>0$.

Let $\lambda = \frac{n+1}{2}$, $N$ a positive integer, $\Sigma$ the set of primes
dividing $N$ and $k$ a
positive integer such that $k > \max\{3, 2\lambda\}$. In
order to define the Eisenstein series we need a Hecke character
$\chi$ of $\mathbb{A}^{\times}$ satisfying
\begin{eqnarray}\label{eqn:character}
\chi_{\infty}(x) &=& \SGN(x)^{k},\\
\chi_{\ell}(a) &=& 1 \quad \text{if $\ell \in \textbf{f}$, $a \in
\inte_{\ell}^{\times}$, and $N \mid (a-1)$}. \nonumber
\end{eqnarray}
Set $D = D[1,N]$ and define functions $\mu$ and $\varepsilon$ on
$\Sp_{2n}(\A)$ by
\begin{eqnarray*}
\mu(x) &=& 0 \quad \text{if $x \notin P(\A)D$},\\
\mu(pw) &=&
 \chi(\det(d_{p}))^{-1}\chi_{\Sigma}(\det(d_{w}))^{-1}\det(d_{p})^{-k}
 \quad \text{if $x=pw \in P(\A)D$,}
\end{eqnarray*}
and
\begin{eqnarray*}
 \varepsilon(x_{\infty}) &=& |j(x_{\infty},i)|^{2}\\
 \varepsilon(x_{\f}) &=& \det(d_{p})^{-2} \quad \text{for} \,\, x=pw
\end{eqnarray*}
where $\displaystyle \chi_{\Sigma} = \prod_{\ell \in \Sigma}\chi_{\ell}$ and
$j(x_{\infty},Z) = \det(c_{x_{\infty}} Z + d_{x_{\infty}})$.

We now have all the ingredients necessary to define the Eisenstein
series we are interested in.  For $x \in \Sp_{2n}(\A)$ and $s \in
\cmplx$, define
\begin{equation*}
E(x,s) = E(x,s; \chi, D) = \sum_{\alpha \in A} \mu(\alpha
x)\varepsilon(\alpha x)^{-s}, \quad A = P\backslash
\Sp_{2n}(\rat).
\end{equation*}
This gives us an Eisenstein series defined on $\Sp_{2n}(\A) \times
\cmplx$, but we will ultimately be interested in an Eisenstein
series $E(Z,s)$ defined on $\h{n} \times \cmplx$.  The Eisenstein series $E(Z,s)$ converges
locally uniformly in $\h{n}$ for $\Real(s) > \lambda$. We associate the Eisenstein series $E(Z,s)$
to $E(x,s)$ as follows.

More generally, let $F_0$ be a function on $\Sp_{2n}(\A)$ such
that
\begin{equation}\label{eqn:blah1}
F_0(\alpha x w) = F_0(x) J(w,i)^{-1} \quad \text{for $\alpha \in
\Sp_{2n}(\rat)$ and $w \in C'$}
\end{equation}
where $C'$ is an open subgroup of $C$ and $J(x,z)$ is defined by
\begin{equation*}
J(x,z) = J_{k,s}(x,z) = j(x,z)^{k}|j(x,z)|^{s}.
\end{equation*}
Our Eisenstein series is such a function. Let $\Gamma' =
\Sp_{2n}(\rat) \cap \Sp_{2n}(\real) C'$ and define a function $F$
on $\h{n}$ by
\begin{equation}\label{eqn:blah2}
F(x(i)) = F_0(x) J(x,i) \quad \text{for $x \in
\Sp_{2n}(\real) C'$}.
\end{equation}
Using the strong approximation theorem ($\Sp_{2n}(\A) =
\Sp_{2n}(\rat) \Sp_{2n}(\real) C'$) we have that $F$ is
well-defined and satisfies
\begin{equation}\label{eqn:blah3}
F(\gamma Z) = F(Z) J(\gamma, Z) \quad \text{for $\gamma \in
\Gamma'$ and $Z \in \h{n}$}.
\end{equation}
Therefore, we have an associated Eisenstein series $E(Z,s)$
defined on $\h{n} \times \cmplx$. The Eisenstein
series $E(Z,s)$ converges locally uniformly in $\h{n}$ for
$\Real(s) > \lambda$. Conversely, given a function
$F$ satisfying Equation \ref{eqn:blah3}, we can define a function
$F_0$ satisfying Equation \ref{eqn:blah1} and Equation
\ref{eqn:blah2} by
\begin{equation*}
F_0(\alpha x) = F(x(i)) J(x,i)^{-1} \quad \text{for $\alpha \in
\Sp_{2n}(\rat)$ and $x \in \Sp_{2n}(\real) C'$}.
\end{equation*}
We will also make use of the fact that if $G =
F\Slash{\gamma^{-1}}$ for $\gamma \in \Gamma_{n}$ with $F$ a
Siegel modular form, then $G_0(x) = F_0(x\gamma_{\textbf{f}})$ and
vice versa.  

\subsection{The Fourier coefficients of $E(Z,s,\chi)$}

We will now focus our attention on the Fourier coefficients of
$E(x,s)$ and in turn $E(Z,s)$.  It turns out that it is
easier to study the Fourier coefficients of a simple translation
of $E(x,s)$ given by
\begin{equation*}
E^{*}(x,s) = E(x \iota_{\textbf{f}}^{-1}, s; \chi, D)
\end{equation*}
where we recall $\iota = \mat{0_{n}}{-1_{n}}{1_{n}}{0_{n}}$ (\cite{shimuraeseries}). Using
the discussion above, we get a corresponding form $E^{*}(Z,s)$.

Let $L=\mathbb{S}^{n}(\rat)\cap \M_{n}(\inte)$, $L' = \{\mathfrak{s} \in
 \mathbb{S}^{n}(\rat):\Tr(\mathfrak{s}L)\subseteq \inte\}$ and $M = N^{-1}L'$.
The Eisenstein series $E^{*}(Z,s)$ has a Fourier expansion
 \begin{equation*}
 E^{*}(Z,s) = \sum_{h\in M} a(h,Y,s)
 e(\Tr(h X))
 \end{equation*}
 for $Z = X+ i Y\in \h{n}$ (\cite{shimuraeseries}).

 \begin{remark}  The Fourier coefficients of
 $E^{*}(Z,s)$ are nonvanishing only when $h$ is totally positive
 definite due to the fact that we have restricted our $k$ to be
 larger then 3 (\cite{shimuraeseries}, Page 460).
 \end{remark}

 We have the following result of Shimura explicitly calculating the Fourier coefficients $a(h,Y,s)$. 

 \begin{proposition}\label{thm:es3} (\cite{shimurabook}, Prop. 18.7, 18.14) For $N\neq 1$,
 \begin{eqnarray*}
 \lefteqn{\hspace*{-1in} a(h,Y,s) = \det(Y)^{-k/2} N^{-n\lambda}
 \det(Y)^{s}\alpha_{N}( ^{t}\!(\overline{Y^{1/2}})h Y^{1/2}; 2s,\chi)}\\ & & \mbox{} \cdot \xi(Y, h,
       s+k/2, s-k/2)
 \end{eqnarray*}
 where $\xi$ is defined by
 \begin{equation*}
 \xi(Y,h;s,t) =
 \int_{\mathbb{S}^{n}(\real)}e(-\Tr(hX))\det(X+iY)^{-s}
 \det(X-iY)^{-t} dX
 \end{equation*}
 with $0<Y\in \mathbb{S}^{n}(\real), \, h\in \mathbb{S}^{n}(\real)$, $s,t \in
 \cmplx$ and $\alpha_{N}$ is a Whittaker integral.  One can consult \cite{shimurabook} for the
 definition of $\alpha_{N}$; it will not be needed here.
 \end{proposition}

 For a Dirichlet character $\psi$, set $\displaystyle \Lambda^{\Sigma}(s,\psi) =
 L^{\Sigma}(2s,\psi)\prod_{j=1}^{[n/2]}L^{\Sigma}(4s-2j,\psi^2)$.
 We normalize $E^{*}(Z,s)$ by multiplying it by
$\pi^{-\frac{n(n+2)}{4}} \Lambda^{\Sigma}(s,\chi)$
 and call this normalized Eisenstein series $D_{E^{*}}(Z,s) =
 D_{E^{*}}(Z,s;k,\chi,N)$.
Consider the Fourier expansion of $D_{E^{*}}(Z,s)$ at $s= \lambda - k/2$:
 \begin{eqnarray*}
 D_{E^{*}}(Z, \lambda-k/2) &=& \sum_{h\in M}
 \pi^{-\frac{n(n+2)}{4}}
 \Lambda^{\Sigma}(\lambda-k/2,\chi) a(h,Y,\lambda-k/2) e(\Tr(hX))\\
    &=& \sum_{h\in M} b(h,Y,\lambda-k/2) e(\Tr(hX)).
 \end{eqnarray*}
 The normalized Eisenstein series $D_{E^{*}}(Z, \lambda -
k/2)$ is in $\mathcal{M}_{k}(\rat^{\text{ab}})$ where
$\rat^{\text{ab}}$ is the maximal abelian extension of $\rat$  (\cite{shimuranearlyholo}, Prop. 4.1).
We show that the coefficients of $D_{E^{*}}(Z, \lambda - k/2)$
actually lie in a finite extension of $\inte_{p}$ for a suitably
chosen prime $p$. Using (\cite{shimuraconfluent}, 4.34K,
4.35IV) we have that
\begin{equation}\label{eqn:xi} \displaystyle
 \xi(Y,h;\lambda, \lambda-k) = \frac{i^{nk}\pi^{\frac{n(n+2)}{4}}
 2^{n(k-1)} \det(Y)^{k-\lambda}}{\mathcal{P}_{n}} \,e(i\,\Tr(hY)),
 \end{equation}
 where
 \begin{equation*}
 \mathcal{P}_{n} = \prod_{j=0}^{\left[\lambda\right]}j! \prod_{j=0}^{\left[\lambda\right]
 -1} \frac{(2j+1)!!}{2^{j+1}}
 \end{equation*}
 and
 \begin{equation*}
  n!! = \left\{ \begin{array}{lll}
            n(n-2)\dots 5 \cdot 3\cdot 1 && n>0, \,\text{odd}\\
            n(n-2) \dots 6\cdot 4\cdot 2 && n>0, \,\text{even}.
            \end{array}
            \right.
 \end{equation*}
 Using that $h$ is totally positive definite we have:
 \begin{proposition}
 \label{thm:es4}
 (\cite{shimurabook}, Prop. 19.2) Set $\chi_{h}$ to be the Hecke
 character corresponding to $\rat(\sqrt{-\det(h)}\,)/\rat$. Then
 \begin{equation*}
 \alpha_{N}(h,s,\chi) =
 \Lambda^{\Sigma}(s,\chi)^{-1}\Lambda^{\Sigma}_{h}(s,\chi)
 \prod_{\ell \in
 \mathcal{C}}f_{h,Y,\ell}(\chi(\ell)|\ell|^{2s})
 \end{equation*}
 where $\mathcal{C}$ is a finite subset of $\f$, the $f_{h,Y,\ell}$
 are polynomials with a constant term of 1 and coefficients in
 $\inte$ independent of $\chi$, and
 \[ \Lambda_{h}^{\Sigma}(s,\chi) = \left\{\begin{array}{ll}
    L^{\Sigma}(2s-n/2,\chi\chi_{h})& n\in 2\inte\\
    1 & \text{otherwise.}\end{array}
    \right.\]
 \end{proposition}

To ease the notation set $\displaystyle \mathcal{F}_{h,Y}(s,\chi)
= \prod_{\ell \in \mathcal{C}}f_{h,Y,\ell}(\chi(\ell)|\ell|^{s})$.
Combining Equation \ref{eqn:xi}, Corollary \ref{thm:es3}, and
Proposition \ref{thm:es4} we have
\begin{equation*}
b(h,Y, \lambda- k/2) =
 \left\{\begin{array}{ll} \frac{i^{nk}2^{n(k-1)} L^{\Sigma}(2\lambda -k-n/2,\chi\chi_{h})
 \mathcal{F}_{Y,h}(2\lambda-k, \chi)}{N^{n\lambda}\mathcal{P}_{n}}\, e(i\,\Tr(hY)) \quad\quad n \in 2\inte\\
 \frac{i^{nk}2^{n(k-1)}
 \mathcal{F}_{Y,h}(2\lambda-k, \chi)}{N^{n\lambda} \mathcal{P}_{n}}\, e(i\,\Tr(hY)) \hspace*{1.15in} \text{otherwise.}
 \end{array} \right.
\end{equation*}

 Let $p$ be an odd prime with $\GCD(p,N)=1$ and $p > 2\lambda -1$. We show
 that the $b(h,Y,\lambda - k/2)$ all lie in
 $\inte_{p}[\chi,i^{nk}]$ where $\inte_{p}[\chi]$ is the extension of $\inte_{p}$ generated
 by the values of $\chi$.  It is clear that $i^{nk}2^{n(k-1)}
 N^{-n\lambda} \in \inte_{p}[\chi,i^{nk}]$ by our choice
 of $p$.  The fact that $p > 2\lambda-1$ and $n \geq 1$ so that $2\lambda -1 \geq \lambda$
 shows that $\mathcal{P}_{n}$
 is in $\inte_{p}$. The fact that we have chosen $k > 2 \lambda$ gives us that
 $2 \lambda -k < 0$.  This in turn shows that $|p|^{2\lambda -k} = p^{k-2\lambda} \in \inte_{p}$.
 Using this fact and that the coefficients of $f_{h,Y,\ell}$ all lie in
 $\inte$, we have
 that $\mathcal{F}_{Y,h}(2\lambda-k,\chi) \in
 \inte_{p}[\chi,i^{nk}]$ for all $h$.  Therefore it remains to
 show that $L^{\Sigma}(2\lambda - k - n/2,\chi \chi_{h}) \in
 \inte_{p}[\chi,i^{nk}]$. We will in fact show that for any
 Dirichlet character $\psi$ of conductor $N$ and any positive
 integer $n$ that $L^{\Sigma}(1-n,\psi) \in \inte_{p}[\psi]$.

 Let $\omega:\inte_{p}^{\times} \rightarrow \mu_{p-1}$ be the usual
Teichmuller character.  One has the existence of a $p$-adic $L$-function $\mathcal{L}_{p}(s,\chi)$
defined on $\{s \in \cmplx_{p} : |s| < (p-1) p^{-1/(p-1)}\}$ such that
\begin{equation*}
\mathcal{L}_{p}(1-n,\psi) = (1 - \psi\omega^{-n}(p)p^{n-1})\,
\frac{B_{n,\psi\omega^{-n}}}{n}
\end{equation*}
for $n\geq 1$ (\cite{washington}, Theorem 5.11).
Using this and the well-known fact that one has
$L(1-n,\psi) =
 -\frac{B_{n,\psi}}{n}$ where $B_{n,\psi}$ is the generalized
 Bernoulli number defined by
 \begin{equation*}
 \sum_{a=1}^{N} \frac{\psi(a)te^{at}}{e^{Nt} - 1} = \sum_{j=0}^{\infty} B_{j,\psi}
 \frac{t^{j}}{j!},
 \end{equation*}
 we can write
\begin{equation*}
L^{\Sigma}(1-n,\psi) = -(1-\psi(p)p^{n-1})^{-1} \prod_{\ell \mid N}
(1-\psi(\ell)\ell^{1-n}) \mathcal{L}_{p}(1-n,\psi\omega^{n}).
\end{equation*}
One can see that $(1-\psi(p)p^{n-1})^{-1} \in \inte_{p}[\psi]$ by
expanding it in a convergent geometric series.  We use the fact
that $\GCD(p,N)=1$ to conclude that $\displaystyle \prod_{\ell
\mid N} (1-\psi(\ell)\ell^{1-n})$ lies in $\inte_{p}[\psi]$. To finish
our proof that $L^{\Sigma}(1-n,\psi) \in \inte_{p}[\psi]$ for all $n
\in \nat$, we note that $\mathcal{L}_{p}(m,\psi)$ is a $p$-adic
integer for all $m$ and all $\psi$ with conductor $N$ such that
$\GCD(p,N)=1$ by (\cite{washington}, Corl. 5.13).  Therefore we
have proven:

\begin{theorem}\label{thm:eiscoeff} Let $n$, $N$, and $k$ be positive integers such
that $k > \max \{3, n+1 \}$. Let $\chi$ be a Dirichlet
character as in Equation \ref{eqn:character}.  Let $p$ be an odd prime
such that $p > n$ and $(p, N) =1$.  Then $D_{E^{*}}(Z,(n+1)/2 -
k/2)$ is in
$\mathcal{M}_{k}(\Gamma_0^{n}(N),\inte_{p}[\chi,i^{nk}])$ for 
\begin{equation*}
\Gamma_0^{n}(N) = \left\{ \gamma \in \Gamma_{n}: c_{\gamma} \equiv 0 (\modu N)\right\}.
\end{equation*}
\end{theorem}

\subsection{Pullbacks and an inner product relation}

In this section we will use the results in the previous section
specialized to the case $n=4$.

We turn our attention to studying the pullback of the Eisenstein
series $E(\mathfrak{Z},s,\chi)$ via maps
\begin{eqnarray*}
\h{2}\times \h{2} &\hookrightarrow& \h{4}\\
(Z,W) &\mapsto& \mat{Z}{0}{0}{W} = \diag[Z,W]
\end{eqnarray*}
and
\begin{eqnarray*}
\Gamma_2 \times \Gamma_2 &\hookrightarrow& \Gamma_4\\
(\alpha, \beta) &\mapsto& \alpha \times \beta = \begin{pmatrix}
a_{\alpha} & 0 & b_{\alpha} & 0 \\ 0 & a_{\beta} & 0 & b_{\beta}
\\ c_{\alpha} & 0 & d_{\alpha} & 0 \\ 0 & c_{\beta} & 0 &
d_{\beta} \end{pmatrix}.
\end{eqnarray*}
These pullbacks have been studied extensively by Shimura
(\cite{shimurazfctns}, \cite{shimurabook}) as well as by Garrett
(\cite{garrett1}, \cite{garrett2}).  In particular, if one has a
Siegel modular form $G$ on $\Gamma_4$ of weight $k$ and level
$N$, then its pullback to $\Gamma_2 \times \Gamma_2$ is a
Siegel modular form in each of the variables $Z$ and $W$ of weight
$k$ and level $N$. We will be interested primarily in the results
found in \cite{shimurazfctns}, particularly the inner product
relation found there.

Let $\sigma_{\f} \in \Sp_8(\rat_{\textbf{f}})$ be defined as
$\sigma_{\f} = (\sigma_{\ell})$ with
\begin{equation*}
\sigma_{\ell} = \left \{ \begin{array}{ccc} I_{8} & & \text{if $\ell \nmid N$} \\
\begin{pmatrix}I_4 & 0_4 \\ \begin{pmatrix} 0_2 & I_2 \\ I_2 &
0_2 \end{pmatrix} & I_4 \end{pmatrix} & & \text{if $\ell \mid N$.}
\end{array}\right.
\end{equation*}
The strong approximation gives an element $\rho \in
\Gamma_4\cap D[1,N]\sigma_{\f}$ such that
\linebreak$N_{\ell}\mid a(\sigma_{\f} \rho^{-1})_{\ell} - I_4$ for
every $\ell \mid N$. In particular, we have that $E\Slash{\rho}$
corresponds to $E(x\sigma_{\f}^{-1})$.

Let $F \in \mathcal{S}_{k}(\Gamma_0^2(N),\real)$ be a Siegel
eigenform.  We specialize a result of
Shimura that gives the inner product of $E\Slash{\rho}$ with such
an $F$. Applying (\cite{shimurazfctns}, Equation 6.17) to our
situation we get
\begin{equation}\label{eqn:shimuraeq1}
\langle D_{E\Slash{\rho}}(\diag[Z,W], (5- k)/2),
(F\Slash{\iota})^{c}(W) \rangle = \pi^{-3} \,\mathcal{A}_{k,N}
L^{\Sigma}_{\st}(5-k,F, \chi) F(Z)
\end{equation}
where $\displaystyle \mathcal{A}_{k,N} =  \frac{(-1)^{k}\,
2^{2k-3} v_{N}}{3 \,[\Gamma_2:\,\Gamma_0^2(N)]}$,
$v_{N} = \pm 1$,
$L^{\Sigma}_{\st}(5-k,F,\chi)$ is the standard zeta function as
defined in Equation \ref{eqn:standardzeta}, and
$(F\Slash{\iota})^{c}$ denotes taking the complex conjugates of
the Fourier coefficients of $F\Slash{\iota}$ where
$F\Slash{\iota}$ is now a Siegel form on
\begin{equation*}
\Gamma^{2,0}(N) = \left\{ \mat{A_2}{B_2}{C_2}{D_2} \in
\Gamma_2| B_2 \equiv 0 (\modu N)\right\}.
\end{equation*}
We can use the $q$-expansion principle for Siegel modular forms
(\cite{chaifaltings}, Prop. 1.5) to conclude that
$F\Slash{\iota}$ has real Fourier coefficients since we chose $F$
to have real Fourier coefficients.  Therefore
$(F\Slash{\iota})^{c}(W)$ in Equation \ref{eqn:shimuraeq1} becomes
$(F\Slash{\iota})(W)$.  Thus we have
\begin{equation*}
\langle D_{E\Slash{\rho (1 \times
\iota_2^{-1})}}(\diag[Z,W], (5-k)/2), F(W) \rangle = \pi^{-3} \,\mathcal{A}_{k,N}
L^{\Sigma}_{\st}(5-k,F, \chi) F(Z).
\end{equation*}

Our next step is to make sure that the Fourier coefficients of
$\mathcal{E}(Z,W)$ are still in some finite extension of
$\inte_{p}$, where
\begin{equation*}
\mathcal{E}(Z,W):= D_{E\Slash{\rho (1 \times
\iota_2^{-1})}}(\diag[Z,W], (5-k)/2).
\end{equation*}
Recall from Theorem \ref{thm:eiscoeff} that $D_{E^{*}}(Z, (5-k)/2)
\in \mathcal{M}_{k}(\Gamma_0^4(N), \inte_{p}[\chi])$. Therefore,
applying the $q$-expansion principle (\cite{chaifaltings},
Prop. 1.5) to $D_{E^{*}}(\diag[Z,W], (5-k)/2)$ slashed by
$\iota_{4}^{-1} \rho (1 \times \iota_{2}^{-1})$, we get that
$D_{E\Slash{\rho} (1 \times \iota_2^{-1})}(\diag[Z,W], (5-k)/2)$
has Fourier coefficients in $\inte_{p}[\chi]$.

Summarizing, we have the following theorem.

\begin{theorem}\label{thm:eisenstein} Let $N>1$ and $k>3$.
For $F \in S_{k}(\Gamma_0^2(N),\real)$ a Hecke eigenform and $p$ a
prime with $p >2$ and $\GCD(p,N)=1$ we have
\begin{equation}\label{eqn:shimurainnerproduct}
\langle \mathcal{E}(Z,W), F(W) \rangle =
\pi^{-3}\,\mathcal{A}_{k,N} L^{\Sigma}_{\st}(5-k,F, \chi) F(Z)
\end{equation}
with $\mathcal{E}(Z,W)$ having Fourier coefficients in
$\inte_{p}[\chi]$.
\end{theorem}

\section{Periods and a certain Hecke operator}\label{sec:periods}

Throughout this section we make the following assumptions.  Let
$k$ be a positive integer with $k \geq 2$. Let $p$ be a
prime so that $p
> 2k-2$.  We let $K$ be a finite extension of $\rat_{p}$ with ring
of integers $\ROI$ and uniformizer $\varpi$. Fix an embedding of $K$ into
$\cmplx$ compatible with the embeddings fixed in Section \ref{sec:notation}. Let $\p$ be the
prime of $\ROI$ lying over $p$.

\subsection{Periods associated to newforms}\label{subsec:periods}

Let $f \in S_{2k-2}(\Gone)$ be a newform with
eigenvalues in $\ROI$. The congruence class of $f$ modulo $p$ is the set of eigenforms
with eigenvalues congruent to those of $f$ modulo $p$. The congruence class of $f$ in
$S_{2k-2}(\Gone)$ determines a maximal ideal $\mathfrak{m}$
of $\T_{\ROI}$ and a residual representation
\begin{equation*}
\rho_{\mathfrak{m}} : \Gal(\overline{\rat}/\rat) \rightarrow
\GL_2(\T_{\ROI}/\mathfrak{m}),
\end{equation*}
so that $\Tr(\rho_{\mathfrak{m}}(\Frob_{\ell})) = T(\ell)$ for all
primes $\ell \neq p$ where $\T_{\ROI}/\mathfrak{m}$ is of
characteristic $p$. This fact is essentially due to Deligne, see
(\cite{ribet1}, Prop. 5.1) for a detailed proof. 

Associated to $f$ is a surjective $\ROI$-algebra map
$\pi_{f}:\T_{\ROI, \mathfrak{m}} \rightarrow \ROI$ given by
$T(\ell) \mapsto a_{f}(\ell)$.  We can view this as a map into
$\cmplx$ as well via the embeddings $\ROI \hookrightarrow K
\hookrightarrow \cmplx$ where the embedding of $K$ into $\cmplx$
was fixed at the beginning of this section.  Let $\wp_{f}$ be the
kernel of $\pi_{f}$.

For $f$ so that $\rho_{\mathfrak{m}}$ is irreducible, one has complex periods $\Omega_{f}^{\pm}$ uniquely determined
up to a $\ROI$-unit as defined in \cite{vatsalcanonical}.  One should
note that while Vatsal restricts to the case of level $N \geq 4$ in \cite{vatsalcanonical}, one
can also define the periods $\Omega_{f}^{\pm}$ for all levels by using the arguments given in (\cite{hidasugaku}, Section 3).
Using these periods we have the following theorem essentially due to Shimura.

\begin{theorem} (\cite{shimuraperiods}, Theorem 1)\label{thm:periods} Let $f \in S_{2k-2}(\Gone, \ROI)$
be a newform. There exist
complex periods $\Omega_{f}^{\pm}$ such that for each integer $m$
with $0 < m < 2k-2$ and every Dirichlet character $\chi$ one has
\begin{equation*}
\frac{L(m,f,\chi)}{\tau(\chi) (2\pi i)^{m}} \in \left\{
\begin{array}{ccc}
\Omega_{f}^{+} \ROI_{\chi} & \text{if} & \chi(-1) = (-1)^{m}\\
 \Omega_{f}^{-}\ROI_{\chi} &\text{if} &\chi(-1)=
(-1)^{m-1}, \end{array} \right.
\end{equation*}
where $\tau(\chi)$ is the Gauss sum of $\chi$ and $\ROI_{\chi}$ is
the extension of $\ROI$ generated by the values of $\chi$.
\end{theorem}

Using the periods $\Omega_{f}^{\pm}$ we make the following conjecture which we 
prove under the additional assumption that $f$ is ordinary at $p$.

\begin{conjecture}\label{conj:heckeop} Let $f=f_1, f_2, \dots f_{r}$ be a basis of
eigenforms for $S_{2k-2}(\Gone)$ with $f$ a newform. Enlarge $\ROI$ if necessary so that the basis
is defined over $\ROI$.  Let $\mathfrak{m}$ be the maximal ideal
in $\T_{\ROI}$ associated to $f$ and assume that representation
$\rho_{\mathfrak{m}}$ is irreducible.  Then there exists a Hecke
operator $t \in \T_{\ROI}$ so that
\begin{equation*}
t f_{i} = \left\{ \begin{array} {ll} u \frac{\langle f, f
\rangle}{\Omega_{f}^{+}\Omega_{f}^{-}} f & \text{if $i = 1$}\\ 0 &
\text{if $i \neq 1$} \end{array} \right.
\end{equation*}
for $u$ a unit in $\ROI$.
\end{conjecture}

\subsection{A certain Hecke operator}

In this section we will establish the validity of Conjecture
\ref{conj:heckeop} in the case that $f$ is ordinary at $p$.

Let $f=f_1, f_2, \dots, f_{r}$ be a basis of
eigenforms for $S_{2k-2}(\Gone)$ as in Conjecture \ref{conj:heckeop}. 
We enlarge $K$ here if
necessary so that this basis is defined over $\ROI$. As with $f$,
there are maps $\pi_{f_{i}}$ for each $i$ as well as kernels
$\wp_{f_{i}}$.

The fact that $f$ is a newform allows us to write
\begin{equation*}
\T_{\ROI, \mathfrak{m}} \otimes_{\ROI} K = K \oplus D
\end{equation*}
for a $K$-algebra $D$ so that $\pi_{f}$ induces the projection of
$\T_{\ROI, \mathfrak{m}}$ onto $K$ (\cite{hidasugaku}).  In this direct sum, $K$
corresponds to the Hecke algebra acting on the eigenspace
generated by $f$ and $D$ corresponds to the Hecke algebra acting
on the space generated by the rest of the $f_{i}$'s. Let $\varrho$
be the projection map of $\T_{\ROI, \mathfrak{m}}$ to $D$.  Set
$I_{f}$ to be the kernel of $\varrho$. Using that our Hecke
algebra is reduced, it is clear from the definition that we have
\begin{equation}\label{eqn:defofif}
I_{f} = \Ann(\wp_{f}) = \bigcap_{i=2}^{r} \wp_{f_{i}}
\end{equation}
where $\Ann(\wp_{f})$ denotes the annihilator of the ideal
$\wp_{f}$.  Since $\T_{\ROI,\mathfrak{m}}$ is reduced, we have
that $\wp_{f} \cap I_{f} = 0$.  Therefore we have that
\begin{equation*}
\T_{\ROI,\mathfrak{m}}/(\wp_{f} \oplus I_{f}) =
\T_{\ROI,\mathfrak{m}}/(\wp_{f}, I_{f}) \xrightarrow{\simeq}
\ROI/\pi_{f}(I_{f})
\end{equation*}
where we use here that
\begin{equation*}
\pi_{f}: \T_{\ROI,\mathfrak{m}}/\wp_{f} \xrightarrow{\simeq} \ROI.
\end{equation*}
Since $\ROI$ is a principal ideal domain, there exists $a \in
\ROI$ so that $\pi_{f}(I_{f}) = a \ROI$.  Therefore we have
\begin{equation}\label{eqn:congmoduleisom}
\ROI/a\ROI \cong \T_{\ROI,\mathfrak{m}}/(\wp_{f} \oplus I_{f}).
\end{equation}

For each prime $\ell$, choose $\alpha_{f}(\ell)$ and
$\beta_{f}(\ell)$ so that $\alpha_{f}(\ell) + \beta_{f}(\ell) =
a_{f}(\ell)$ and \linebreak$\alpha_{f}(\ell) \beta_{f}(\ell) =
\ell^{2k-3}$. Set
\begin{equation*}
D(s,\pi_{f}) = \prod_{\ell} \left( ( 1- \alpha_{f}(\ell)^2 \ell^{-s})(1-
\alpha_{f}(\ell)\beta_{f}(\ell) \ell^{-s}) (1 - \beta_{f}(\ell)^2
\ell^{-s}) \right)^{-1}.
\end{equation*}
Shimura has shown this Euler product converges if the real part of
$s$ is sufficiently large and can be extended to a meromorphic
function on the entire complex plane that is holomorphic except
for possible simple poles at $s=2k-2$ and $2k-3$
(\cite{shimuraproc}, Theorem 1).  The values of
$D(2k-2,\pi_{f})/U(\pi_{f})$ are in $\ROI$ (\cite{hidasugaku},
Page 86) where
\begin{equation*}
U(\pi_{f}) = \frac{(2 \pi)^{2k-1}\,
\Omega_{f}^{+}\Omega_{f}^{-}}{(2k-3)!}.
\end{equation*}
Following Hida we define $\varepsilon \in K$ by
\begin{equation}\label{eqn:a1}
a = \frac{D(2k-2,\pi_{f})}{\varepsilon \cdot U(\pi_{f})}
\end{equation}
where $a$ is given by Equation \ref{eqn:congmoduleisom}.

\begin{theorem}\label{thm:hidaepsilon} (\cite{hidasugaku}, Theorem 2.5)  Let $f \in
S_{2k-2}(\Gone, \ROI)$ be a newform. Let $\p$ be the
prime of $\ROI$ over $p$. If $f$ is ordinary at $\p$, then
$\varepsilon$ is a unit in $\ROI$.
\end{theorem}

Combining (\cite{hidacongruences}, Theorem 5.1) and
(\cite{shimuraarithmetic}, 8.2.17) we have
\begin{equation*}
D(2k-2, \pi_{f}) = \frac{2^{4k-4}\, \pi^{2k-1}}{(2k-3)!} \, \langle f, f \rangle.
\end{equation*}
Inserting this expression
for $D(2k-2,\pi_{f})$ into Equation \ref{eqn:a1} and simplifying
we obtain
\begin{equation*}
a = \frac{2^{2k-3}}{\varepsilon \cdot
\Omega_{f}^{+}\Omega_{f}^{-}} \langle f, f \rangle.
 \end{equation*}
Combining Equations \ref{eqn:defofif} and \ref{eqn:congmoduleisom}
we can write
\begin{equation}
\T_{\ROI, \mathfrak{m}}/(\wp_{f} \oplus \bigcap_{i=2}^{r}
\wp_{f_{i}}) \cong \ROI/a \ROI
\end{equation}
where
\begin{equation}\label{eqn:a}
a = \frac{2^{2k-3}}{\varepsilon \cdot
 \Omega_{f}^{+}\Omega_{f}^{-}} \langle f, f \rangle.
\end{equation}
Since $\T_{\ROI, \mathfrak{m}}/\wp_{f} \cong \ROI$, there exists a
$t \in I_{f}$ that maps to $a$ under the above isomorphism. Thus
we have that
\begin{equation*}
t f_{i} = \left \{ \begin{array}{ccc} a f & & \text{if $i=1$}\\
0 & & \text{if $2 \leq i \leq r$.} \end{array} \right.
\end{equation*}
This is the Hecke operator we seek. Using the fact that
\begin{equation*}
\T_{\ROI} \cong \prod \T_{\ROI,\mathfrak{m}}
\end{equation*}
where the product is over the maximal ideals of $\T_{\ROI}$, we
can view $\T_{\ROI,\mathfrak{m}}$ as a subring of $\T_{\ROI}$.
Therefore we have the following theorem.

\begin{theorem}\label{thm:heckeoperator} Let $f= f_1, f_2, \dots, f_{r}$ be a
basis of eigenforms of $S_{2k-2}(\Gone, \ROI)$ with $k
>2$. Suppose that the representation
$\rho_{\mathfrak{m}}$ associated to $f$ is irreducible and $f$ is
ordinary at $\p$. There exists a Hecke operator $t \in \T_{\ROI}$
such that $tf = a f$ and $tf_{i} = 0$ for $i \geq 2$ where $a$ is
as in Equation \ref{eqn:a} with $\varepsilon$ a unit in $\ROI$.
\end{theorem}

\section{The congruence}\label{chap:congruence}

In this section we combine the results of the
previous sections to produce a congruence between the
Saito-Kurokawa lift $\Ff$ and a cuspidal Siegel eigenform $G$ which is not
a Saito-Kurokawa lift.
We fix $k>3$ throughout this section.

\subsection{Congruent to a Siegel modular form}\label{sec:congmodform}

Let $f \in S_{2k-2}(\Gone)$ be a newform and $\Ff$ the Saito-Kurokawa lift 
as constructed in Section \ref{sec:sk}.  
Recall that $\mathcal{E}(Z,W)$ is a Siegel modular form of weight
$k$ and level $N$ in each variable.  Before we go any further we 
need to replace $\mathcal{E}(Z,W)$ with a form
of level $1$. The reason for this will be clear shortly as we will
need to apply a Hecke operator that is of level $1$. We do this by
taking the trace. Set
\begin{equation*}
\tilde{\mathcal{E}}(Z,W) = \sum_{\gamma \times \delta \in
\Gtwo/\Gamma_0^2(N)\times \Gtwo/\Gamma_0^2(N)}
\mathcal{E}(Z,W)\Slash{(\gamma\times \delta)}.
\end{equation*}
It is clear that $\tilde{\mathcal{E}}(Z,W)$ is now a Siegel modular
form on $\Gtwo \times \Gtwo$.  The Fourier
coefficients are seen to still be in $\inte_{p}[\chi]$ by applying
the $q$-expansion principle for Siegel modular forms
(\cite{chaifaltings}, Prop. 1.5). 

Let $F_0 = \Ff, F_1, \dots F_{r}$ be a basis of eigenforms for the
Hecke operators $T(\ell)$ ($\ell \neq p$) of
$\mathcal{M}_{k}(\Gtwo)$ so that $F_{i}$ is orthogonal
to $\Ff$ for $1\leq i \leq r$. We enlarge
$\ROI$ here if necessary so that \\
1. $\ROI$ contains the values of $\chi$\\
2. the eigenforms $F_{i}$ are all defined over $\ROI$\\
3. the newforms $f_{i}$ defined in Conjecture \ref{conj:heckeop}
are defined over $\ROI$.\\
Following Shimura, we write
\begin{equation}\label{eqn:eisensteinsum}
\tilde{\mathcal{E}}(Z,W) = \sum_{i,j} c_{i,j} F_{i}(Z)
F_{j}(W)
\end{equation}
with $c_{i,j} \in \cmplx$ (\cite{shimurazfctns}, Eq. 7.7).

\begin{lemma}\label{thm:removingf} Equation
\ref{eqn:eisensteinsum} can be written in the form
\begin{equation*}
\tilde{\mathcal{E}}(Z,W) = c_{0,0}\Ff(Z) \Ff(W) +
\sum_{\small{\begin{array}{c} 0\leq i \leq r\\0<j\leq
r\end{array}}} c_{i,j}F_{i}(Z) F_{j}(W).
\end{equation*}
\end{lemma}

\begin{proof} Recall Shimura's inner product formula as given in
Equation \ref{eqn:shimurainnerproduct}:
\begin{equation*}
\langle \mathcal{E}(Z,W), \Ff(W) \rangle_{\Gamma^2_0(N)} =
\pi^{-3}\,\mathcal{A}_{k,N} L_{\st}^{\Sigma}(5-k,\Ff, \chi) \Ff(Z)
\end{equation*}
and observe that
\begin{equation*}
\langle \mathcal{E}(Z,W), \Ff(W) \rangle_{\Gamma^2_0(N)} = \langle
\tilde{\mathcal{E}}(Z,W), \Ff(W) \rangle_{\Gtwo}
\end{equation*}
by the way we defined the inner product. Note that we insert the
``$\Gamma^2_0(N)$'' and ``$\Gtwo$'' here merely to make
explicit which group the inner product is defined on.
 On the other hand, if we take the inner product of the right
hand side of Equation \ref{eqn:eisensteinsum} with $\Ff(W)$ we get
\begin{equation*}
\langle \tilde{\mathcal{E}}(Z,W), \Ff(W) \rangle = \sum_{0 \leq i \leq
r} c_{i,0} \langle \Ff, \Ff \rangle F_{i}(Z).
\end{equation*}
Equating the two we get
\begin{equation*}
\pi^{-3}\,\mathcal{A}_{k,N} L_{\st}^{\Sigma}(5-k,\Ff, \chi) \Ff(Z)
= \sum_{0 \leq i \leq r} c_{i,0} \langle \Ff, \Ff \rangle
F_{i}(Z).
\end{equation*}
Since the $F_{i}$ form a basis, it must be the case that $c_{i,0}
= 0$ unless $i=0$, which gives the result.
\end{proof}

Our goal is to show
that we can write $c_{0,0}$ as a product of a unit in $\ROI$ and
$\frac{1}{\varpi^{m}}$ for some $m \geq 1$. Once we have shown we
can do this, it will be straightforward to move from this to the
congruence we desire.

Using Lemma \ref{thm:removingf} and Equation
\ref{eqn:shimurainnerproduct} we write
\begin{equation}
c_{0,0} \langle \Ff,\Ff \rangle \Ff(Z) =
\pi^{-3}\,\mathcal{A}_{k,N} L_{\st}^{\Sigma}(5-k, \Ff, \chi) \Ff(Z).
\end{equation}
Equating the coefficient of $\Ff(Z)$ on each side and solving for
$c_{0,0}$ gives us
\begin{equation}\label{eqn:c001}
c_{0,0} = \frac{\mathcal{A}_{k,N}
L_{\st}^{\Sigma}(5-k,\Ff,\chi)}{\pi^3 \,\langle \Ff, \Ff \rangle}.
\end{equation}
Combining Theorems \ref{thm:innerproducts} and \ref{thm:standardzetafactorization} with Equation \ref{eqn:c001} we have
\begin{equation}
c_{0,0} = \mathcal{B}_{k,N}\frac{|D|^{k-3/2}\,
L(k-1,f,\chi_{D})L^{\Sigma}(3-k,\chi) L^{\Sigma}(1,f,\chi)
L^{\Sigma}(2,f,\chi)}{\pi^2\, |c_{g}(|D|)|^2\, L(k,f)\langle f, f
\rangle}
\end{equation}
with
\begin{equation}
\mathcal{B}_{k,N} =
\frac{(-1)^{k} 2^{2k+2}\, 3 \,v_{N}}{
(k-1) [\Gtwo:\Gamma_0^2(N)]}.
\end{equation}

The main obstacle at this point to studying the $\varpi$-valuation
of $c_{0,0}$ is the possibility that the congruence we produce would be to
a Saito-Kurokawa lift. Fortunately, we
can apply the results of Section \ref{sec:periods} to remove
this possibility.  We will do this by applying a Hecke operator
$t_{S}$ to Equation \ref{eqn:eisensteinsum}.

Assume that Conjecture \ref{conj:heckeop} is satisfied. Recall
that we showed this is the case if $f$ is
ordinary at $\p$. We have a Hecke operator $t \in \T_{\ROI}$ that
acts on $f$ via the eigenvalue $\displaystyle  u \,\frac{\langle
f, f \rangle}{ \Omega_{f}^{+}\Omega_{f}^{-}}$ for $u$ a unit in
$\ROI$ and kills $f_{i}$ for all other $f_{i}$ in a basis of
newforms for $S_{2k-2}(\Gone,\ROI)$. Using that the
Saito-Kurokawa correspondence is Hecke-equivariant,
we have associated to $t$ a Hecke operator $t_{S} \in \T_{S,\ROI}$
so that
\begin{equation}\label{eqn:siegelhecke}
t_{S}\cdot F_{f_{i}} = \left\{ \begin{array}{ccc}  u
\,\frac{\langle f, f \rangle}{ \Omega_{f}^{+}\Omega_{f}^{-}}\,\Ff
& & \text{for $f_{i} = f$}\\ 0 & & \text{for $f_{i} \neq
f$.}\end{array}\right.
\end{equation}
Applying $t_{S}$ to Equation \ref{eqn:eisensteinsum} as a modular form in $W$ we
obtain
\begin{equation}
t_{S}\tilde{\mathcal{E}}(Z,W) = c_{0,0}' \Ff(Z)
\Ff(W) + \sum_{\small{\begin{array}{c} 0 \leq i \leq r \\ 0 < j
\leq r \end{array}}} c_{i,j}F_{i}(Z) t_{S}F_{j}(W)
\end{equation}
with
\begin{equation}\label{eqn:cvalue1}
c_{0,0}' = u \,\frac{\langle f, f \rangle}{
\Omega_{f}^{+}\Omega_{f}^{-}}\cdot c_{0,0} =
\mathcal{C}_{k,N}\, \frac{|D|^{k-3/2}\,
L(k-1,f,\chi_{D})L^{\Sigma}(3-k,\chi) L^{\Sigma}(1,f,\chi)
L^{\Sigma}(2,f,\chi)}{\pi^2\, |c_{g}(|D|)|^2\,L(k,f) \Omega_{f}^{+}
\Omega_{f}^{-}}
\end{equation}
where
\begin{equation}
\mathcal{C}_{k,N} = u \cdot \mathcal{B}_{k,N}.
\end{equation}
Note that we have killed any $F_{j}$ that is a Saito-Kurokawa lift. 

Our next step is to normalize the $L$-values in Equation
\ref{eqn:cvalue1} so as to obtain algebraic values.  Theorem
\ref{thm:periods} showed that if we
divide $L(m,f,\chi)$ by $\tau(\chi) (2\pi i )^{m}
\Omega_{f}^{\pm}$ we get a value in $\ROI$ where we choose
$\Omega_{f}^{+}$ if $\chi(-1) = (-1)^{m}$ and choose
$\Omega_{f}^{-}$ if $\chi(-1) = (-1)^{m-1}$. It is easy to see
that if $\Omega_{f}^{+}$ is associated to $L(1,f,\chi)$, then
$\Omega_{f}^{-}$ is associated to $L(2,f,\chi)$ and vice versa.
Therefore we have
\begin{equation*}
\frac{L(1,f,\chi) L(2,f,\chi)}{\Omega_{f}^{+} \Omega_{f}^{-}} =
\tau(\chi)^2 (2 \pi i)^3 L_{\alg}(1,f,\chi) L_{\alg}(2,f,\chi).
\end{equation*}
In particular, we have
\begin{equation*}
\frac{L^{\Sigma}(1,f,\chi)L^{\Sigma}(2,f,\chi)}{\Omega_{f}^{+}\Omega_{f}^{-}}
= \frac{\tau(\chi)^2 (2 \pi i)^3
L_{\alg}(1,f,\chi)
L_{\alg}(2,f,\chi)}{L_{\Sigma}(1,f,\chi)L_{\Sigma}(2,f,\chi)}.
\end{equation*}

Next we turn our attention to the ratio $\displaystyle
\frac{L(k-1,f,\chi_{D})}{L(k,f)}$. Since $L(k,f)$ has no
character, we see that we associate
$\Omega_{f}^{+}$ to $L(k,f)$ if $k$ is even and $\Omega_{f}^{-}$ if $k$ is odd.
We need to associate the same period to $L(k-1,f,\chi_{D})$. The
way to accomplish this is to choose $D$ so that $\chi_{D}(-1) =
-1$.  Therefore we have
\begin{equation*}
\frac{L(k-1,f,\chi_{D})}{L(k,f)} =
\frac{\tau(\chi_{D})L_{\alg}(k-1,f,\chi_{D})}{(2 \pi i)
L_{\alg}(k,f)}.
\end{equation*}
Also recall that in Section \ref{sec:eseries} we showed that
$L^{\Sigma}(3-k,\chi) \in \inte_{p}[\chi]$ for $\GCD(p,N)=1$.

Gathering these results together we have:
\begin{equation}
c_{0,0}'
= \mathcal{D}_{k,N,\chi,D} \mathcal{L}(k,f,D,\chi)
\end{equation}
where
\begin{equation*}
\mathcal{L}(k,f,D,\chi) = \frac{L^{\Sigma}(3-k,\chi)
L_{\alg}(k-1,f,\chi_{D})L_{\alg}(1,f,\chi)L_{\alg}(2,f,\chi)}{L_{\alg}(k,f)}
\end{equation*}
and
\begin{equation*}
\mathcal{D}:=\mathcal{D}_{k,N,\chi,D} = \frac{(-1)^{k+1} \, 2^{2k+4}\, 3\, |D|^{k} \tau(\chi_{D})
\tau(\chi)^2}{(k-1) \, [\Gtwo:\Gamma_0^2(N)] |D|^{3/2}\, |c_{g}(|D|)|^2 \,L_{\Sigma}(1,f,\chi)L_{\Sigma}(2,f,\chi)}.
\end{equation*}
Everything in these equations is now algebraic, so it comes down to
studying the $\varpi$-divisibility of each of the terms.  We would
like to show that $\varpi^{m}$ divides the denominator for some $m
\geq 1$ but not the numerator. Note that as long as everything in
the denominator is a $\varpi$-integer, we do not have to worry
about anything written in the denominator contributing a
``$\varpi$'' to the numerator.

We first deal with $\mathcal{D}$.  We know from
Conjecture \ref{conj:heckeop} that $u$ is a unit of $\ROI$ so long
as $\rho_{\mathfrak{m}}$ is irreducible. Under this assumption we
need not worry about $u$. Choosing $p$ relatively
prime to $D$ takes care of the $D$'s that appear.  We also can see
that $\varpi \nmid \tau(\chi)$ and $\varpi \nmid \tau(\chi_{D})$.
For instance, suppose $\varpi \mid \tau(\chi)$. Then this would
imply that $\varpi \mid (\tau(\chi)\overline{\tau(\chi)})^2 =
(\sqrt{N})^2=N$, a contradiction and similarly for $\tau(\chi_{D})$.
Next we need to deal with $\displaystyle
\frac{1}{L_{\Sigma}(1,f,\chi)L_{\Sigma}(2,f,\chi)}$.  Observe
that we can write
\begin{eqnarray*}
\frac{1}{L_{(\ell)}(1,f,\chi)} &=&
\frac{1}{(1-\lambda_{f}(\ell)\ell^{-1} + \ell^{2k-5})} \\
&=& \frac{\ell}{(\ell -\lambda_{f}(\ell) + \ell^{2k-4})}.
\end{eqnarray*}
Since $\ell \mid N$ and $\GCD(p,N)=1$, we have that $p \nmid
\ell$.  It is also clear now that the denominator is in $\ROI$.
Since we can do this for each $\ell \mid N$ and the same argument
follows for $L_{(\ell)}(2,f,\chi)$, we see that $\displaystyle
\frac{1}{L_{\Sigma}(1,f,\chi)L_{\Sigma}(2,f,\chi)}$ cannot
possibly contribute any $\varpi$'s to the numerator.  Recall that
$|c_{g_{f}}(|D|)|^2 \in \ROI$.  Therefore, so long as we choose our $p>2k-2$ and $p$ 
relatively prime to $[\Gtwo:\Gamma_0^2(N)]$, we have that 
$\mathcal{D}$ cannot contribute any
$\varpi$'s to the numerator.

The term $\mathcal{L}(k,f,D,\chi)$ is where the divisibility
assumption enters into our calculations. We assume here that for
some integer $m \geq 1$ we have $\varpi^{m} \mid L_{\alg}(k,f)$
and that if $\varpi^{n}
\parallel
L^{\Sigma}(3-k,\chi)L_{\alg}(k-1,f,\chi_{D})L_{\alg}(1,f,\chi)L_{\alg}(2,f,\chi)$
then $n<m$ so that we end up with a $\varpi$ in the denominator of
$c'_{0,0}$.

Under these assumptions we can write
\begin{equation}\label{eqn:seedofg}
t_{S}\tilde{\mathcal{E}}(Z,W) =
\frac{A}{\varpi^{m-n}}\Ff(Z)\Ff(W) + \sum_{\small{
\begin{array}{c}0 \leq i \leq r \\0<j\leq r \end{array}}} c_{i,j}F_{i}(Z) t_{S}F_{j}(W)
\end{equation}
for some $\varpi$-unit $A$. Recall that Corollary
\ref{thm:skfourier} gave that $\Ff$ has Fourier coefficients in
$\ROI$ and that we can find a $T_0$ so that $\varpi \nmid
A_{\Ff}(T_0)$. This allows us to immediately conclude that we must
have some $c_{i,j} \neq 0$ for at least one of $i, j \neq 0$.
Otherwise we would have $t_{S}\tilde{\mathcal{E}}(Z,W) =
\frac{A}{\varpi^{m-n}} \Ff(Z) \Ff(W)$ and using the
integrality of the Fourier coefficients of $t_{S}\tilde{\mathcal{E}}(Z,W)$ 
we would get $\Ff(Z) \Ff(W)
\equiv 0 (\modu \varpi^{m-n})$, a contradiction.

Recall that by Corollary \ref{thm:skfourier} there exists a $T_0$ so that 
$A_{\Ff}(T_0)$ is in $\ROI^{\times}$.  Expand each side of 
Equation \ref{eqn:seedofg} in terms of $Z$, reduce modulo $\varpi$ and equate the
$T_0^{\text{th}}$ Fourier coefficients.  Using the $\ROI$-integrality of the Fourier coefficients of $
t_{S}\tilde{\mathcal{E}}(Z,W)$ we
obtain:
\begin{equation*}
A_{\Ff}(T_0) \Ff(W) \equiv - \frac{\varpi^{m-n}}{A}\sum_{\small{
\begin{array}{c}0 \leq i \leq r \\0<j\leq r \end{array}}} c_{i,j}
A_{F_{i}}(T_0) t_{S}F_{j}(W) (\modu \varpi^{m-n}),
\end{equation*}
i.e., we have a congruence $\Ff \equiv G (\modu \varpi^{m-n})$ for
$G \in \mathcal{M}_{k}(\Gone)$ where
\begin{equation}\label{eqn:formulag}
G(W)= - \frac{\varpi^{m-n}}{A \cdot A_{\Ff}(T_0)} \sum_{\small{
\begin{array}{c}0 \leq i \leq r \\0<j\leq r \end{array}}} c_{i,j}
A_{F_{i}}(T_0) t_{S}F_{j}(W).
\end{equation}
Since the Hecke operator $t_{S}$ killed all $F_{j}$ ($0<j\leq r$) that came from Saito-Kurokawa lifts,
we have that $G$ is a sum of forms that are not Saito-Kurokawa lifts.

Momentarily we will show how $G$ can be used to produce a non-Saito-Kurokawa
cuspidal eigenform with eigenvalues that are congruent to the eigenvalues
of $\Ff$, but before we do we gather our
results into the following theorem.

\begin{theorem}\label{thm:congruentmodular} Let $k>3$ be an integer and $p$
a prime so that $p >2k-2$.  Let $f \in
S_{2k-2}(\Gone,\ROI)$ be a newform with real Fourier
coefficients and $\Ff$ the Saito-Kurokawa lift of $f$. Suppose
that $\rho_{\mathfrak{m}}$ is irreducible and that Conjecture
\ref{conj:heckeop} is satisfied. If there exists an integer $N>1$, a fundamental discriminant $D$ so that $(-1)^{k-1}D>0$,
$\chi_{D}(-1)=-1$, $p\nmid ND[\Gtwo:\Gamma_0^2(N)]$, and a Dirichlet character $\chi$ of conductor $N$ so that
\begin{equation*}
\varpi^{m} \mid L_{\alg}(k,f)
\end{equation*}
with $m \geq 1$ and
\begin{equation*}
\varpi^{n} \parallel
L^{\Sigma}(3-k,\chi)L_{\alg}(k-1,f,\chi_{D})L_{\alg}(1,f,\chi)L_{\alg}(2,f,\chi)
\end{equation*}
with $n < m$, then there exists $G \in
\mathcal{M}_{k}(\Gtwo)$ that is a sum of eigenforms that are not Saito-Kurokawa lifts so
that
\begin{equation*}
\Ff \equiv G (\modu \varpi^{m-n}).
\end{equation*}
\end{theorem}

\subsection{Congruent to a non-Saito-Kurokawa cuspidal eigenform}

In this section we will show how given a congruence
\begin{equation}\label{eqn:firstcongruence}
\Ff \equiv G(\modu \varpi^{m})
\end{equation}
for $m \geq 1$ as in Theorem \ref{thm:congruentmodular}, we can
find a non-Maass cuspidal eigenform that has the same eigenvalues as $\Ff$
modulo $\varpi$.

\begin{notation} If $F_1$ and $F_2$ have eigenvalues that are congruent modulo
$\varpi$, we will write
\begin{equation*}
F_1 \equiv_{\ev} F_2 (\modu \varpi)
\end{equation*}
where the $\ev$ stands for the congruence being a congruence of
eigenvalues.
\end{notation}

We begin by showing that given a congruence as in Theorem
\ref{thm:congruentmodular}, there must be a non-Saito-Kurokawa eigenform $F$ so that
$F \equiv_{\ev} \Ff (\modu \varpi)$. Once we have shown this, we
will show that we can obtain an eigenvalue congruence to
a cusp form. Applying the first result again we obtain our final
goal of an eigenvalue congruence between $\Ff$ and a cuspidal
eigenform that is not a Saito-Kurokawa lift.

\begin{lemma}\label{thm:congeigenform} Let $G \in \mathcal{M}_{k}(\Gtwo)$ be
as in Equation \ref{eqn:formulag} so that we have the congruence
$G \equiv \Ff (\modu \varpi)$. Then there exists an eigenform $F$
so that $F$ is not a Saito-Kurokawa lift and $\Ff \equiv_{\ev} F
(\modu \varpi)$.
\end{lemma}

\begin{proof} As in Equation \ref{eqn:formulag} write $G = \sum c_{i} F_{i}$ with each
$F_{i}$ an eigenform and $c_{i} \in \ROI$. It
is clear from the construction of $G$ that $F_{i} \neq \Ff$ and
$F_{i}$ is not a Saito-Kurokawa lift for all $i$. Recall that we
have the decomposition
\begin{equation*}
\T_{S, \ROI} \cong \prod \T_{S,\ROI, \mathfrak{m}}
\end{equation*}
where the $\mathfrak{m}$ are maximal ideals of $\T_{S,\ROI}$
containing $\varpi$. Let $\mathfrak{m}_{\Ff}$ be the maximal ideal
corresponding to $\Ff$. There is a Hecke operator $t \in
\T_{S,\ROI}$ so that $t \Ff = \Ff$ and $t F = 0$ for any eigenform
$F$ that does not correspond to $\mathfrak{m}_{\Ff}$, i.e., if $F
\not \equiv_{\ev} \Ff (\modu \varpi)$. If $F_{i} \not \equiv_{\ev}
\Ff (\modu \varpi)$ for every $i$ then applying $t$ to the
congruence $G \equiv \Ff (\modu \varpi)$ would then yield $\Ff
\equiv 0 (\modu \varpi)$, a contradiction to the fact that $A_{\Ff}(T_0) \in \ROI^{\times}$. Thus there must
be an $i$ so that $\Ff \equiv_{\ev} F_{i} (\modu \varpi)$.
\end{proof}

We now show that we actually have an eigenvalue congruence to a
cusp form.  Before we prove this fact, we briefly recall the Siegel operator
$\Phi$. The Siegel operator is defined by
\begin{equation*}
\Phi(F(\tau)) = \lim_{\lambda \rightarrow \infty} F\left(\mat{\tau}{0}{0}{i\lambda}\right)
\end{equation*}
where $\tau \in \h{1}$.  In terms of Fourier coefficients we have
\begin{equation*}
\Phi(F(\tau)) = \sum_{n\geq 0} a_{F}\left(\mat{n}{0}{0}{0}\right) e^{2 \pi i n \tau}.
\end{equation*}
From this expression it is clear that if $F$ has Fourier coefficients in $\ROI$, so does $\Phi(F)$.  
We note the following facts about the Siegel operator which can all be found in \cite{freitag}:
\begin{enumerate}
\item Given a Siegel modular form $F \in \mathcal{M}_{k}(\Gtwo)$, one has $\Phi(F)\in M_{k}(\Gone)$.  \\
\item If $\Phi(F) = 0$, then $F$ is a cusp form.\\
\item If $F$ is an eigenform of the operator $T_{S}(\ell)$, then $\Phi(F)$ is an eigenform
of $T(\ell)$.\\
\item \label{item4} We have the following formula:
\begin{equation*}\label{eq:siegeloperator}
\Phi(T_{S}(\ell)F) = (1-\ell^{2-k})T(\ell)\Phi(F).
\end{equation*}
\end{enumerate}

Let $\Ff \equiv_{\ev} F (\modu \varpi)$ with $F$ the non
Saito-Kurokawa eigenform constructed in Lemma
\ref{thm:congeigenform}. Suppose that $F$ is not a cusp form so
that $\Phi(F) \neq 0$. Let $g = \Phi(F)$.
We denote the $n^{\text{th}}$ eigenvalue of $g$ as
$\lambda_{g}(n)$. Let $\ell$ be a prime so that $\ell \neq p$.
Note that since $F$ has eigenvalues in $\ROI$ and (\ref{item4}) gives that $\lambda_{F}(\ell) =
(1-\ell^{2-k})\lambda_{g}(\ell)$, we must have
$(1-\ell^{2-k})\lambda_{g}(\ell) \in \ROI$. Applying (\ref{item4}) again gives
\begin{equation*}
\Phi(T_{S}(\ell) F) 
    = (1-\ell^{2-k})\lambda_{g}(\ell)g.
\end{equation*}
On the other hand, the (\ref{item4}) and the congruence give us that
\begin{eqnarray*}
\Phi(T_{S}(\ell)F) &\equiv_{\ev}& \Phi(T_{S}(\ell)\Ff) (\modu \varpi)\\
    &=& \Phi(\lambda_{\Ff}(\ell) \Ff)\\
    &=& \lambda_{\Ff}(\ell) g\\
    &=& (\ell^{k-1} + \ell^{k-2}+\lambda_{f}(\ell))\, g.
\end{eqnarray*}
Thus we have that
\begin{equation}\label{eqn:galrep1}
(\ell^{k-1} + \ell^{k-2}+\lambda_{f}(\ell))  \equiv
(1-\ell^{2-k})\lambda_{g}(\ell) (\modu \varpi).
\end{equation}

Denote the Galois representation associated to $f$ by $\rho_{f}$
and similarly for $g$. Denote the residual representations after
reducing modulo $\varpi$ by $\overline{\rho}_{f}$ and
$\overline{\rho}_{g}$. Equation \ref{eqn:galrep1} and the
Tchebotarov Density Theorem show that we have the following
equivalence of 4-dimensional Galois representations
\begin{equation*}
\begin{pmatrix}\omega^{k-1} & &  \\ & \omega^{k-2}
& \\  &  & \overline{\rho}_{f} \end{pmatrix} =
\mat{\overline{\rho}_{g}}{}{}{\omega^{2-k}\overline{\rho}_{g}}.
\end{equation*}
It is clear from this that $\overline{\rho}_{f}$
must be reducible. However, we assumed before that this was not
the case. This contradiction shows that $\Phi(F) =0$.  We have
proved the following theorem.

\begin{theorem}\label{thm:thecongruence} Let $k>3$ be an integer and $p$
a prime so that $p >2k-2$.  Let $f \in
S_{2k-2}(\Gone,\ROI)$ be a newform with real Fourier
coefficients and $\Ff$ the Saito-Kurokawa lift of $f$. Suppose
that $\rho_{\mathfrak{m}}$ is irreducible and that Conjecture
\ref{conj:heckeop} is satisfied. If there exists an integer $N>1$, a fundamental discriminant $D$ so that $(-1)^{k-1}D>0$,
$\chi_{D}(-1)=-1$, $p\nmid ND[\Gtwo:\Gamma_0^2(N)]$, and a Dirichlet character $\chi$ of conductor $N$ so that
\begin{equation*}
\varpi^{m} \mid L_{\alg}(k,f)
\end{equation*}
with $m \geq 1$ and
\begin{equation*}
\varpi^{n} \parallel
L^{\Sigma}(3-k,\chi)L_{\alg}(k-1,f,\chi_{D})L_{\alg}(1,f,\chi)L_{\alg}(2,f,\chi)
\end{equation*}
with $n < m$, then there exists an eigenform $G \in
\mathcal{S}_{k}(\Gtwo)$ that is not a
Saito-Kurokawa lift so that
\begin{equation*}
\Ff \equiv_{\ev} G (\modu \varpi).
\end{equation*}
\end{theorem}

\section{Generalities on Selmer groups}\label{sec:selmer1}

In this section we define the relevant Selmer group following
Bloch and Kato \cite{blochkato} and Diamond, Flach, and Guo
\cite{diamondflachguo}. We also collect various results
that are not easily located in existing references. We conclude the section by stating a
version of the Bloch-Kato conjecture for modular forms.

For a field $K$ and a topological $\Gal(\overline{K}/K)$-module
$M$, we write $\cohom^1(K, M)$ for
$\cohom^1_{\cont}(\Gal(\overline{K}/K), M)$ to ease notation,
where ``cont'' indicates continuous cocycles. We write $D_{\ell}$
to denote a decomposition group at $\ell$ and $I_{\ell}$ to denote an
inertia group at $\ell$. We identify $D_{\ell}$
with $\Gal(\overline{\rat}_{\ell}/\rat_{\ell})$.

Let $E$ be a finite extension of $\rat_{p}$, $\ROI$ the ring of
integers of $E$, and $\varpi$ a uniformizer.  Let $V$ be a
$p$-adic Galois representation defined over $E$. Let $T \subseteq
V$ be a $\absgal$-stable $\ROI$-lattice. Set $W = V/T$.  For $n
\geq 1$, put
\begin{equation*}
W_{n} = W [\varpi^{n}]= \left\{ x \in W : \varpi^{n} x = 0
\right\} \cong T/\varpi^{n} T.
\end{equation*}

In the following section we will construct non-zero cohomology
classes in $\cohom^1(\rat, W_1)$ and we would like to know that
they remain non-zero when we map them into
$\cohom^1(\rat, W)$ under the natural map.

\begin{lemma} If $T/\varpi T$ is irreducible as an $\left(\ROI/\varpi \ROI\right)[\absgal]$-module
then $\cohom^1(\rat, W_1)$ injects into $\cohom^1(\rat,W)$.
\end{lemma}

\begin{proof} Consider the exact sequence

\begin{figure}[h!]
\centerline{\xymatrix{ 0 \ar[r]& W_1 \ar[r]& W \ar[r]^{\cdot
\varpi}& W \ar[r] & 0.}}
\end{figure}

\noindent This short exact sequence gives rise to the long exact
sequence of cohomology groups

\begin{figure}[h!]
\centerline{\xymatrix{ 0 \ar[r]& \cohom^0(\rat,W_1) \ar[r]&
\cohom^0(\rat,W) \ar[r]^{\cdot \varpi}& \cohom^0(\rat,W) \ar[r] &\\ &
\cohom^1(\rat,W_1) \ar[r]^{\psi} & \cohom^1(\rat, W) \ar[r]&
\cdots .}}
\end{figure}

\noindent We show that $\psi$ is injective. Recalling that
$\cohom^0(G,M) = M^{G}$, it is clear that $\cohom^0(\rat, W_1) =
0$ since we have assumed that $T/\varpi T$ is irreducible.
Since $W$ is torsion, $\cohom^0(\rat, W)$ is necessarily torsion
as well.  If $\cohom^0(\rat,W)$ contains a non-zero element,
multiplying by a suitable $\varpi^{m}$ makes it a non-zero element
in $W_1$.  This would give us a non-zero element in
$\cohom^0(\rat,W_1)$, a contradiction. Thus, we obtain that
$\cohom^0(\rat,W) =0$ and so $\psi$ is an injection.
\end{proof}

We also have that $\cohom^1(\rat_{\ell}, W_1)$ injects into
$\cohom^1(\rat_{\ell}, W)$ when $T/\varpi T$ is irreducible as an
$\left(\ROI/\varpi \ROI\right)[D_{\ell}]$-module by an analogous
argument.

We write $\B_{\crys}$ to denote Fontaine's ring of $p$-adic
periods as defined in \cite{fontaine1}.  For a $p$-adic representation $V$,
set
\begin{equation*}
D_{\crys} = (V \otimes_{\rat_{p}} \B_{\crys})^{D_{p}}.
\end{equation*}

\begin{definition} A $p$-adic representation $V$ is called {\it
crystalline} if
\begin{equation*}
\dim_{\rat_{p}} V = \dim_{\rat_{p}} D_{\crys}.
\end{equation*}
\end{definition}

\begin{definition} A crystalline representation $V$ is called {\it
short} if the following hold\\
1. $\Fil^0 D_{\crys} = D_{\crys}$ and $\Fil^{p}D_{\crys} =0$,\\
2. if $V'$ is a nonzero quotient of $V$, then $V'
\otimes_{\rat_{p}} \rat_{p}(p-1)$ is ramified\\
where $\Fil^{i} D_{\crys}$ is a decreasing filtration of $D_{\crys}$ as given in \cite{diamondflachguo}.
\end{definition}


Following Bloch-Kato (\cite{blochkato}), we define
spaces $\cohom^1_{f}(\rat_{\ell}, V)$ by
\begin{equation*}
\cohom_{f}^1(\rat_{\ell},V) = \left\{ \begin{array}{lll}
\cohom^1_{\ur}(\rat_{\ell}, V) & & \ell \neq p, \infty\\
\ker(\cohom^1(\rat_{p},V) \rightarrow \cohom^1(\rat_{p}, V \otimes
\B_{\crys})) & & \ell = p
\end{array} \right.
\end{equation*}
where
\begin{equation*}
\cohom^1_{\ur}(\rat_{\ell}, M) = \ker (\cohom^1(\rat_{\ell}, M)
\rightarrow \cohom^1(I_{\ell}, M))
\end{equation*}
for any $D_{\ell}$-module $M$. The Bloch-Kato groups
$\cohom_{f}^1(\rat_{\ell}, W)$ are defined by
\begin{equation*}
\cohom_{f}^1(\rat_{\ell}, W) = \im (\cohom_{f}^1(\rat_{\ell}, V)
\rightarrow \cohom^1(\rat_{\ell}, W)).
\end{equation*}

One should note here that the $f$ appearing in these definitions
has nothing to do with the elliptic modular form $f$ we have been
working with and is merely standard notation in the literature
(standing for ``finite part''.)

\begin{lemma}\label{thm:unramified} If $V$ is unramified at $\ell$, then
\begin{equation*}
\cohom_{f}^1(\rat_{\ell}, W) = \cohom_{\ur}^1(\rat_{\ell}, W).
\end{equation*}
\end{lemma}

\begin{proof} We need only show that
$\cohom^1_{\ur}(\rat_{\ell},V)$ surjects onto
$\cohom^1_{\ur}(\rat_{\ell},W)$. The short exact sequence\\

\begin{figure}[h!]
\centerline{\xymatrix{ 0 \ar[r]& T \ar[r]& V \ar[r] & W \ar[r]& 0
}}
\end{figure}

\noindent gives rise to the long exact sequence in cohomology\\

\begin{figure}[h!]
\centerline{\xymatrix{ 0 \ar[r]& \cohom^0(\F_{\ell}, T) \ar[r]&
\cohom^0(\F_{\ell}, V) \ar[r]& \cohom^0(\F_{\ell}, W) \ar[r]&
\cohom^1(\F_{\ell},T) \ar[r]&\\ & \cohom^1(\F_{\ell}, V) \ar[r]^{\psi}&
\cohom^1(\F_{\ell},W) \ar[r]& \cohom^2(\F_{\ell},T) \ar[r]&
\cdots}}
\end{figure}

\noindent where we identify $\Gal(\overline{\F}_{\ell}/\F_{\ell})$
with $D_{\ell}/I_{\ell}$.  Since
$\Gal(\overline{\F}_{\ell}/\F_{\ell}) \cong \hat{\inte}$ and
$\hat{\inte}$ has cohomological dimension 1
(\cite{serregaloiscohomology}, Chap. 5), we have
that $\cohom^2(\F_{\ell}, T) = 0$, i.e., $\psi$ is a surjection.
Observing that for any $D_{\ell}$-module $M$ we have a natural
isomorphism
\begin{equation*}
\cohom^1_{\ur}(\rat_{\ell}, M) \cong \cohom^1(\F_{\ell},
M^{I_{\ell}})
\end{equation*}
and using the fact that $V$ is assumed to be unramified at $\ell$
and so $T$ is unramified at $\ell$ as well, we obtain the result.
\end{proof}

\begin{remark}\label{rem:cocyle} Let $R$ be a ring and let $M$ and $N$ be
$R$-modules.  Recall that an ($R$-linear) extension of $M$ by $N$
is a short
exact sequence of $R$-modules

\begin{figure}[h!]
\centerline{\xymatrix{0 \ar[r]& N \ar[r]& X \ar[r]& M \ar[r]& 0.}}
\end{figure}

\noindent There is a bijection between $\Ext_{R}^1(M,N)$ and the
set of equivalence classes of extensions of $M$ by $N$. Let $\alpha \in
\cohom^1(\rat_{\ell}, V)$.  It is known that
$\cohom^1(\rat_{\ell}, V) \cong \Ext_{E[D_{\ell}]}^1(E,V)$ (\cite{jacobson}, Theorem 6.12). Therefore we have that $\alpha$
corresponds to an extension $X$ of $E$ by $V$:

\begin{figure}[h!]
\centerline{\xymatrix{ 0 \ar[r]& V \ar[r]& X \ar[r]& E \ar[r] &
0.}}
\end{figure}

\noindent For $\ell \neq p$, one has that $X$ is an unramified
representation if and only if \linebreak$\alpha \in
\cohom^{1}_{\ur}(\rat_{\ell}, V)$.  If $\ell = p$, then $X$ is a
crystalline representation if and only if $\alpha \in
H^1_{f}(\rat_{p}, V)$.
\end{remark}

We are now in a position to define the Selmer group of interest to
us.

\begin{definition}\label{def:selmar}  Let $W$ and $\cohom_{f}^1(\rat_{\ell}, W)$ be defined as above.
The {\it Selmer group} of $W$ is given by
\begin{equation*}
\cohom^1_{f}(\rat, W) = \ker \left( \cohom^1(\rat,W) \rightarrow
\bigoplus_{\ell} \frac{\cohom^1(\rat_{\ell},
W)}{\cohom^1_{f}(\rat_{\ell}, W)}\right),
\end{equation*}
i.e., it consists of the cocycles $c \in \cohom^1(\rat, W)$ that
when restricted to $D_{\ell}$ lie in $\cohom^1_{f}(\rat_{\ell},
W)$ for each $\ell$.
\end{definition}

Lemma \ref{thm:unramified} allows us to identify
$\cohom_{f}^1(\rat_{\ell}, W)$ with $\cohom_{\ur}^1(\rat_{\ell},
W)$ for $\ell \neq p$.  Define $\cohom_{f}^1(\rat_{\ell}, W_{n}) =
\cohom_{\ur}^1(\rat_{\ell}, W_{n})$ for $\ell \neq p$. At the
prime $p$, we define $\cohom_{f}^1(\rat_{p}, W_{n}) \subseteq
\cohom^1(\rat_{p}, W_{n})$ to be the subset of classes of
extensions of $D_{p}$-modules

\begin{figure}[h!]
\centerline{\xymatrix{ 0 \ar[r]& W_{n} \ar[r]& X \ar[r]&
\ROI/\varpi^{n}\ROI \ar[r] & 0}}
\end{figure}

\noindent so that $X$ is in the essential image of $\mathbb{V}$
where $\mathbb{V}$ is the functor defined in Section 1.1 of
\cite{diamondflachguo}.  We will not define the functor here; we
will be content with stating the relevant properties that we will
need.  This essential image is stable under direct sums,
subobjects, and quotients (\cite{diamondflachguo}, Section 2.1).
This gives that $\cohom_{f}^1(\rat_{p}, W_{n})$ is an
$\ROI$-submodule of $\cohom^1(\rat_{p},W_{n})$.  We also have that
$\cohom_{f}^1(\rat_{p}, W_{n})$ is the preimage of
$\cohom_{f}^1(\rat_{p},W_{n+1})$ under the natural map
$\cohom^1(\rat_{p}, W_{n}) \rightarrow \cohom^1(\rat_{p},
W_{n+1})$. For our purposes, it will be enough to note the
following fact.

\begin{lemma} (\cite{diamondflachguo}, Page
670)\label{thm:diamondflachguo2} If $V$ is a short crystalline
representation at $p$, $T$ a $D_{p}$-stable lattice, and $X$ a
subquotient of $T/\varpi^{n}T$ that gives an extension of
$D_{p}$-modules as above then the class of this extension is in
$\cohom_{f}^1(\rat_{p}, W_{n})$.
\end{lemma}


We have a natural map $\phi_{n}: \cohom^1(\rat_{p}, W_{n})
\rightarrow \cohom^1(\rat_{p}, W)$.  On the level of extensions
this map is given by pushout via the map $\varpi^{-n}T/T
\rightarrow V/T$, pullback via the map $\ROI \rightarrow
\ROI/\varpi^{n}\ROI$, and the isomorphism $\cohom^1(\rat_{p}, W)
\cong \Ext^1_{\ROI[D_{p}]}(\ROI, V/T)$. In the next section we
will be interested in the situation where we have a non-zero
cocycle $h \in \cohom^1(\rat, W_1)$ that restricts to be in
$\cohom^1_{f}(\rat_{\ell}, W_1)$.  We want to be able to conclude
that this gives a non-zero cocycle in $\cohom^1(\rat, W)$ that
restricts to be in $\cohom^1_{f}(\rat_{\ell}, W)$.  We saw above
that $\cohom^1(\rat, W_1)$ injects into $\cohom^1(\rat, W)$, so it
only remains to show that the restriction is in
$\cohom^1_{f}(\rat_{\ell}, W)$. This is accomplished via the
following proposition.

\begin{proposition}(\cite{diamondflachguo}, Prop. 2.2)\label{thm:diamondflachguo} The
natural isomorphism
\begin{equation*}
\varinjlim_{n} \cohom^1(\rat_{\ell}, W_{n}) \cong
\cohom^1(\rat_{\ell}, W)
\end{equation*}
induces isomorphisms
\begin{equation*}
\varinjlim_{n} \cohom_{\ur}^1(\rat_{\ell}, W_{n}) \cong
\cohom_{\ur}^1(\rat_{\ell},W)
\end{equation*}
and
\begin{equation*}
\varinjlim_{n} \cohom_{f}^1(\rat_{p}, W_{n}) \cong
\cohom_{f}^1(\rat_{p}, W).
\end{equation*}
\end{proposition}

\noindent This proposition shows that the map $\phi_{n}$ gives a
map from $\cohom^1_{f}(\rat_{p}, W_{n})$ to
$\cohom^1_{f}(\rat_{p}, W)$.  We summarize with the following proposition.

\begin{proposition}\label{thm:selmersummary}  Let $h$ be a non-zero cocycle in
$\cohom^1(\rat, W_{1})$ and assume that $T/\varpi T$ is
irreducible.  If $h|_{D_{\ell}} \in \cohom^1_{f}(\rat_{\ell},
W_{1})$ is non-zero, then $h|_{D_{\ell}}$ gives a non-zero
$\varpi$-torsion element of $\cohom^1_{f}(\rat_{\ell}, W)$.  If
$h|_{D_{\ell}} \in \cohom^1_{f}(\rat_{\ell}, W_{1})$ for every
prime $\ell$, then $h$ is a non-zero $\varpi$-torsion element of
$\cohom_{f}^1(\rat, W)$.
\end{proposition}

We conclude this section with a brief discussion of the Bloch-Kato
conjecture for modular forms.  The reader interested in more
details or a more general framework should consult
\cite{blochkato} where the conjecture is referred to as the
``Tamagawa number conjecture''.

For each prime $p$ let $V_{p}:=V_{f,\p}$ be the $p$-adic Galois
representation arising from a newform $f$ of weight $2k-2$,
$T_{p}:=T_{f,\p}$ a $\absgal$-stable lattice, and $W_{p}:=
W_{f,\p}= V_{p}/T_{p}$. The $W_{p}$ here should not be confused with our
use of $W_{n}$ earlier.  Denote the $j^{\text{th}}$ Tate twist of
$W_{p}$ by $W_{p}(j)$. Let $\pi_{*}$ be the natural map
$\cohom^1(\rat,V_{p}(j)) \rightarrow \cohom^1(\rat,W_{p}(j))$ used
to define the groups $\cohom^1_{f}(\rat_{\ell},W_{p}(j))$. We
define the Tate-Shafarevich group to be
\begin{equation}
\Sha(j) = \bigoplus_{\ell}
\cohom^1_{f}(\rat,W_{\ell}(j))/\pi_{*}\cohom_{f}^1(\rat,
V_{\ell}(j)).
\end{equation}
Define the set $\Gamma_{\rat}(j)$ by
\begin{equation*}
\Gamma_{\rat}(j) = \bigoplus_{\ell} \cohom^0(\rat, W_{\ell}(j)).
\end{equation*}
One should think of these as the analogue of the rational torsion
points on an elliptic curve. Accordingly, the set $\Gamma_{\rat}(j)$ is often referred
to as the ``global points''.

\begin{conjecture}(Bloch-Kato)\label{thm:bkconjecture} With the
notation as above, one has
\begin{equation}\label{eqn:blochkatoequation}
L(k,f) = \frac{\left(\prod_{\ell} c_{\ell}(k)\right)
\vol_{\infty}(k) \# \Sha(1-k)}{\# \Gamma_{\rat}(k) \#
\Gamma_{\rat}(k-2)}
\end{equation}
where $c_{p}(j)$ are ``Tamagawa factors'' and $\vol_{\infty}(k)$
is a certain real period. See \cite{deligne} for a careful
treatment of $\vol_{\infty}(k)$.
\end{conjecture}

\begin{remark}\label{remark:blochkato} 1.  It is known that away from the central critical
value the Selmer group is finite (\cite{kato}, Theorem 14.2).
Therefore we can identify the $\varpi$-part of the Selmer group
with the $\varpi$-part of the Tate-Shafarevich group.\\
2.  If $T_{\ell}/\varpi T_{\ell}$ is irreducible, then
$\cohom^0(\rat, W_{\ell}(j)) = 0$. \\
3. The Tamagawa factors are integers.  See (\cite{blochkato}) for definitions and discussion.\\
4. The real period $\vol_{\infty}(k)$ is $\pi^{k}
\Omega_{f}^{\pm}$ up to $p$-adic unit with the $\pm$ depending on
the parity of $k$ (\cite{dummigan}).
\end{remark}

In the next section we will prove that if $\varpi \mid
L_{\alg}(k,f)$, then $p \mid \# \cohom^1_{f}(\rat,W_{f,\p}(1-k))$.
Using Remark \ref{remark:blochkato} this divisibility gives
evidence for the Bloch-Kato conjecture as stated.  In particular, we
will have that if a prime $\varpi$ divides the left hand side of Equation 
\ref{eqn:blochkatoequation}, then it divides the right hand side as well.

\section{Galois Arguments}

In this section we will combine the results of the previous
sections to imply a divisibility result on the Selmer group
$\cohom_{f}^1(\rat, W_{f,\p}(1-k))$.  Note that in this section entries of matrices denoted by $*$'s can
be anything and are assumed to be of the appropriate size.
Similarly, a $1$ as a matrix entry is assumed to be of the
appropriate size.  A blank space in a matrix is assumed to be 0.  We begin by 
stating two theorems that are fundamental to the results in this section.

\begin{theorem}\label{thm:siegelgalois} (\cite{skinnerurban}, Theorem 3.1.3)
Let $F \in \mathcal{S}_{k}(\Gamma^2_0(M))$ be an eigenform, $K_{F}$
the number field generated by the Hecke eigenvalues of $F$, and
$\p$ a prime of $K_{F}$ over $p$.  There exists a finite extension
$E$ of the completion $K_{F,\p}$ of $K_{F}$ at $\p$ and a
continuous semi-simple Galois representation
\begin{equation*}
\rho_{F,\p}: \absgal \rightarrow \GL_4(E)
\end{equation*}
unramified at all primes $\ell \nmid pM$ so that for all $\ell
\nmid pM$, we have
\begin{equation*}
\det(X\cdot I - \rho_{F,\p}(\Frob_{\ell})) = L_{\spin,(\ell)}(X)
\end{equation*}
(we are using arithmetic Frobenius here as opposed to geometric
which is more prevalent in the literature.)
\end{theorem}

\begin{theorem}\label{thm:crystallineshort} (\cite{faltings2}, \cite{urbandecomp}) Let $F$ be as in Theorem
\ref{thm:siegelgalois}. The restriction of $\rho_{F,\p}$ to the decomposition group
$D_{p}$ is crystalline at $p$. In addition if $p
> 2k-2$ then $\rho_{F,\p}$ is short.
\end{theorem}

Recall that for a Saito-Kurokawa lift one has a decomposition of
the Spinor $L$-function: for $\Ff$ we have
\begin{equation*}
L_{\spin}(s,\Ff) = \zeta(s-k+1)\zeta(s-k+2)L(s,f).
\end{equation*}
This decomposition gives us that the Galois representation
$\rho_{\Ff,\p}$ has a very simple form.  In particular, using that
$\rho_{\Ff,\p}$ is semi-simple and applying the Brauer-Nesbitt
theorem we have that
\begin{equation*}
\rho_{\Ff,\p} = \begin{pmatrix} \varepsilon^{k-2} & & \\ &
\rho_{f,\p}&\\ & & \varepsilon^{k-1} \end{pmatrix}
\end{equation*}
where $\varepsilon$ is the $p$-adic cyclotomic character.

Under the conditions of Theorem \ref{thm:thecongruence}, we have a non-Saito-Kurokawa
cuspidal Siegel eigenform $G$ such that $G \equiv_{\ev} \Ff (\modu
\varpi)$. This gives a congruence between the Hecke polynomials of the Spinor
$L$-functions of $\Ff$ and $G$ as well. Let $\p$ be a prime of a sufficiently large
finite extension $E/\rat_{p}$ so that $\ROI_{E}$ contains the
$\ROI$ needed for the congruence and so that $\rho_{\Ff,\p}$ and
$\rho_{G,\p}$ are defined over $\ROI_{E}$.  We set $\ROI =
\ROI_{E}$ and let $\varpi$ be a uniformizer of $\ROI$.  Applying the 
Brauer-Nesbitt theorem we obtain that the semi-simplification of $\overline{\rho}_{G,\p}$ 
is given by
\begin{equation*}
\overline{\rho}_{G,\p}^{\semi} = \overline{\rho}_{\Ff,\p} =
\begin{pmatrix} \omega^{k-2} & & \\ & \overline{\rho}_{f,\p}&\\ & &
\omega^{k-1} \end{pmatrix}
\end{equation*}
where we use $\omega$ to denote the reduction of the cyclotomic
character $\varepsilon$ modulo $\varpi$.  The goal now is to use
this information on the semi-simplification of
$\overline{\rho}_{G,\p}$ to deduce the form of
$\overline{\rho}_{G,\p}$.

Our first step is to show that there is a $\absgal$-stable lattice
$T$ so that the reduction of $\rho_{G,\p}$ is of the form
\begin{equation*}
\overline{\rho}_{G,\p} = \begin{pmatrix} \omega^{k-2} & *_1 & *_2
\\ *_3 & \overline{\rho}_{f,\p} & *_4 \\ & &
\omega^{k-1}\end{pmatrix}
\end{equation*}
where either $*_1$ or $*_3$ is zero. We proceed by brute force,
working our way backwards from the definition of the semi-simplification.
We begin by noting some conjugation formulas that will be
important. Expanding on the notation used in \cite{ribet1}, write
\begin{equation*}
P_1 = \begin{pmatrix} \varpi & & & \\ & 1 & & \\ & & 1 & \\
& & & 1\end{pmatrix},
\end{equation*}
\begin{equation*}
P_2 = \begin{pmatrix} 1 & & & \\ & \varpi & & \\ & & \varpi & \\&
& & 1
\end{pmatrix},
\end{equation*}
and
\begin{equation*}
P_3 = \begin{pmatrix} 1& & \\ & 1&& \\ & & 1 &  \\ & & &\varpi
\end{pmatrix}.
\end{equation*}
We have the following conjugation formulas
\begin{equation}\label{eqn:conjugation1}
P_1 \begin{pmatrix} a_{1,1} &a_{1,2}& a_{1,3} & a_{1,4}\\ \varpi
a_{2,1}
&a_{2,2}& a_{2,3} & a_{2,4}\\ \varpi a_{3,1} &a_{3,2}& a_{3,3} & a_{3,4}\\
\varpi a_{4,1} &a_{4,2}& a_{4,3} & a_{4,4}\end{pmatrix} P_1^{-1} =
\begin{pmatrix}a_{1,1} &\varpi a_{1,2}& \varpi a_{1,3} & \varpi a_{1,4}\\a_{2,1}
&a_{2,2}& a_{2,3} & a_{2,4}\\a_{3,1} &a_{3,2}& a_{3,3} & a_{3,4}\\
a_{4,1} &a_{4,2}& a_{4,3} & a_{4,4}\end{pmatrix},
\end{equation}
\begin{equation}\label{eqn:conjugation2}
P_2 \begin{pmatrix} a_{1,1} & \varpi a_{1,2}& \varpi a_{1,3} & a_{1,4}\\
a_{2,1}
&a_{2,2}& a_{2,3} & a_{2,4}\\  a_{3,1} &a_{3,2}& a_{3,3} & a_{3,4}\\
 a_{4,1} & \varpi a_{4,2}& \varpi a_{4,3} & a_{4,4}\end{pmatrix} P_2^{-1} =
\begin{pmatrix}a_{1,1} &  a_{1,2}&  a_{1,3} &  a_{1,4}\\ \varpi a_{2,1}
&a_{2,2}& a_{2,3} & \varpi a_{2,4}\\\varpi a_{3,1} &a_{3,2}& a_{3,3} & \varpi a_{3,4}\\
a_{4,1} &a_{4,2}& a_{4,3} & a_{4,4}\end{pmatrix},
\end{equation}
and
\begin{equation}\label{eqn:conjugation3}
P_3 \begin{pmatrix} a_{1,1} &a_{1,2}& a_{1,3} & \varpi
a_{1,4}\\a_{2,1}
&a_{2,2}& a_{2,3} & \varpi a_{2,4}\\a_{3,1} &a_{3,2}& a_{3,3} & \varpi a_{3,4}\\
a_{4,1} &a_{4,2}& a_{4,3} & a_{4,4}\end{pmatrix} P_3^{-1} =
\begin{pmatrix}a_{1,1} &a_{1,2}& a_{1,3} & a_{1,4}\\a_{2,1}
&a_{2,2}& a_{2,3} & a_{2,4}\\a_{3,1} &a_{3,2}& a_{3,3} & a_{3,4}\\
\varpi a_{4,1} & \varpi a_{4,2}& \varpi a_{4,3} &
a_{4,4}\end{pmatrix}.
\end{equation}

The definition of semi-simplification gives us vector spaces 
\begin{equation*}
V := V_{G,\p} = V_0 \supset V_1 \supset V_2
\supset V_3 =0
\end{equation*}
with each of $V_0/V_1$, $V_1/V_2$, and $V_2$
irreducible components of $\overline{\rho}_{G,\p}$.  Since we know
$\overline{\rho}^{\semi}_{G,\p}$ explicitly, we can say that
$V_0/V_1$, $V_1/V_2$, and $V_2$ consist of two 1-dimensional
spaces and one 2-dimensional space, corresponding to
$\omega^{k-1}, \omega^{k-2}$ and $\overline{\rho}_{f,\p}$.  The
difficulty is that we do not know which $V_{i}/V_{i+1}$
corresponds to which of $\omega^{k-1}$, $\omega^{k-2}$, and
$\overline{\rho}_{f,\p}$.  We handle this by considering all possible
situations and seeing what this implies for the form of
$\overline{\rho}_{G,\p}$.  We split this into several cases.\\

\noindent{\bf Case 1:} $\dim V_2 = 1 = \dim V_0/V_1$, $\dim
V_1/V_2 = 2$.\\

\noindent{\bf Case 2:} $\dim V_2 = 2$, $\dim V_0/V_1 = \dim
V_1/V_2 = 1$.\\

\noindent{\bf Case 3:} $\dim V_2 = \dim V_1/V_2 =1$, $\dim V_0/V_1
= 2$.\\

Each of these cases can be analyzed via the conjugation formulas given above. We illustrate this 
with Case 2.  This case corresponds to the situation where we have either
\begin{equation*}
\overline{\rho}_{G,\p} = \begin{pmatrix} \overline{\rho}_{f,\p} &
* & * \\ & \omega^{k-2} & * \\ & & \omega^{k-1} \end{pmatrix}
\end{equation*}
or
\begin{equation*}
\overline{\rho}_{G,\p} = \begin{pmatrix}\overline{\rho}_{f,\p} &
* & * \\ & \omega^{k-1} & * \\ & & \omega^{k-2} \end{pmatrix}.
\end{equation*}

\noindent The first of these is handled by observing
\begin{equation*}
\begin{pmatrix} & 1 & \\ 1 & & \\ & & 1 \end{pmatrix} \begin{pmatrix} \overline{\rho}_{f,\p} &
* & * \\ & \omega^{k-2} & * \\ & & \omega^{k-1} \end{pmatrix} \begin{pmatrix} & 1 & \\ 1 & & \\ & & 1
\end{pmatrix}^{-1} = \begin{pmatrix} \omega^{k-2} & & * \\ * &
\overline{\rho}_{f,\p} & * \\ & & \omega^{k-1} \end{pmatrix}.
\end{equation*}

\noindent The second is handled similarly:
\begin{equation*}
\begin{pmatrix} & & 1 \\ 1 & & \\ & 1 &
\end{pmatrix}\begin{pmatrix}\overline{\rho}_{f,\p} &
* & * \\ & \omega^{k-1} & * \\ & & \omega^{k-2} \end{pmatrix} \begin{pmatrix} & & 1 \\ 1 & & \\ & 1 &
\end{pmatrix}^{-1} = \begin{pmatrix} \omega^{k-2} & & \\ * & \overline{\rho}_{f,\p} & *
\\ * & & \omega^{k-1} \end{pmatrix} .
\end{equation*}
Next we change bases as in Equation \ref{eqn:conjugation2} and
then as in Equation \ref{eqn:conjugation3} to obtain
\begin{equation*}
\overline{\rho}_{G,\p} = \begin{pmatrix} \omega^{k-2} & * & * \\ &
\overline{\rho}_{f,\p} & * \\ & & \omega^{k-1} \end{pmatrix}.
\end{equation*}

Therefore we have that there is a lattice so that we have
\begin{equation*}
\overline{\rho}_{G,\p} = \begin{pmatrix} \omega^{k-2} & *_1 & *_2
\\ *_3 & \overline{\rho}_{f,\p} & *_4 \\ & &
\omega^{k-1}\end{pmatrix}
\end{equation*}
where either $*_1$ or $*_3$ is zero.

Now that we have the matrix in the appropriate form, we would like
to further limit the possibilities.  We begin with the following
proposition.

\begin{proposition}\label{thm:ribettypethm} Let $\rho_{G,\p}$ be such that it does not have a sub-quotient
of dimension 1 and $\overline{\rho}_{G,\p}^{\semi} = \omega^{k-2}
\oplus \overline{\rho}_{f,\p}\oplus \omega^{k-1}$. Then there
exists a $\absgal$-stable $\ROI$-lattice in $V_{G}$ having an
$\ROI$-basis such that the corresponding representation
$\rho=\rho_{G,\p}: \absgal \rightarrow \GL_4(\ROI)$ has reduction
of the form
\begin{equation}\label{eqn:reduction1}
\overline{\rho}_{G,\p} = \begin{pmatrix} \omega^{k-2} & *_1 & *_2
\\ *_3 & \overline{\rho}_{f,\p} & *_4 \\ & &
\omega^{k-1}\end{pmatrix}
\end{equation}
and such that there is no matrix of the form
\begin{equation}\label{eqn:u}
U = \begin{pmatrix} 1 & & & n_1\\& 1 & & n_2\\& & 1& n_3 \\& & & 1
\end{pmatrix} \in \GL_4(\ROI)
\end{equation}
such that $\rho' = U \rho U^{-1}$ has reduction of type
(\ref{eqn:reduction1}) with $*_2 = *_4 = 0$.
\end{proposition}

\begin{proof} Fix a $\absgal$-stable lattice and an $\ROI$-basis
giving rise to a representation $\rho_{0}$ of type
(\ref{eqn:reduction1}). Suppose there exists a $U_0$ as in
(\ref{eqn:u}).  Inductively we define a converging sequence
of matrices $M_{i}$ so that $M_{i} \rho_{0} M_{i}^{-1}$ is a 
representation into $\GL_4(\ROI)$ with reduction of the form
(\ref{eqn:reduction1}).  Set $M_1 = U_0$.  By assumption we have that
$M_1 \rho_{0} M_1^{-1}$ is of the required form.  Define $M_{i+1}$ inductively 
by $M_{i+1} = P_{3}^{-i}U_0P_{3}^{i}M_{i}$. 
We have that
\begin{equation*}
M_{i+1} = \begin{pmatrix} 1 & 0& 0& n_1 \sum_{n=1}^{i} \varpi^{n} \\
			0&1&0& n_2 \sum_{n=1}^{i} \varpi^{n}\\
			0&0&1& n_3 \sum_{n=1}^{i} \varpi^{n}\\
			0&0&0&1
\end{pmatrix}.
\end{equation*}
From this it is clear that $M_{i}$ converges to some $M_{\infty} \in \GL_{4}(\ROI)$
of the form
\begin{equation*}
M_{\infty} = \begin{pmatrix} 1 & 0& 0& t_1 \\
			0&1&0& t_2\\
			0&0&1& t_3\\
			0&0&0&1
\end{pmatrix}
\end{equation*}
where $\displaystyle t_{j} = n_{j} \lim_{i \rightarrow \infty} \sum_{n=1}^{i} \varpi^{n}$.
Suppose we have that $M_{i} \rho_{0} M_{i}^{-1}$ is of the required form. 
Using the defintion of $M_{i+1}$ we have that $M_{i+1} \rho_0
M_{i+1}^{-1}$ is of the form that the first three entries of the
rightmost column are all divisible by $\varpi^{i}$ since $P_3^{i} M_{i+1} \rho_0
M_{i+1}^{-1} P_{3}^{-i}$ has entries in $\ROI$. Thus,
$\rho_{\infty} = M_{\infty} \rho_0 M_{\infty}^{-1}$ is such that
the first three entries of the rightmost column are zero.  This
gives a 1-dimensional subquotient of $\rho_{G,\p}$, a
contradiction.  Thus no such $U_0$ can exist.
\end{proof}

In light of this proposition our next step is to show that
$\rho_{G,\p}$ does not have a sub-quotient of dimension 1
as in Theorem \ref{thm:thecongruence}.  There are three
possibilities for how $\rho_{G,\p}$ could split up with a
sub-quotient of dimension 1.  It could have a sub-quotient of
dimension 3 and of dimension 1, a 2-dimensional sub-quotient and
two 1-dimensional ones, or four 1-dimensional sub-quotients.  The
case of a 3 dimensional sub-quotient cannot occur, see
(\cite{urbanduke}, Page 512) or (\cite{skinnerurban}, Proof of
Theorem 3.2.1).  The case of splitting into four 1-dimensional
sub-quotients is not possible either.  Indeed, if $\rho_{G,\p} =
\chi_1 \oplus \chi_2 \oplus \chi_3 \oplus \chi_4$ for characters
$\chi_{i}$, then $\overline{\rho}_{G,\p}$ splits into four
1-dimensional sub-quotients as well but this gives a contradiction
as we know $\overline{\rho}_{f,\p}$ is not completely reducible
(\cite{ribet1}, Prop. 2.1).

The last case to worry about is if $\rho_{G,\p}$ splits into a
2-dimensional sub-quotient and two 1-dimensional sub-quotients. In
this case $G$ must be a CAP form (\cite{skinnerurban}, Proof of
Theorem 3.2.1) induced from the Siegel parabolic. However, the 
results of \cite{shapiro} imply that $G$ must then be a Saito-Kurokawa lift,
a contradiction.

Summarizing to this point, we now have that there exists a
$\absgal$-stable lattice $T_{G,\p}$ so that the reduction
$\overline{\rho}_{G,\p}$ is of the form
\begin{equation*}
\overline{\rho}_{G,\p} = \begin{pmatrix} \omega^{k-2} & *_1 & *_2
\\ *_3 & \overline{\rho}_{f,\p} & *_4 \\ & & \omega^{k-1}
\end{pmatrix}
\end{equation*}
where $*_1$ or $*_3$ is zero and so that $\overline{\rho}_{G,\p}$
is not equivalent to a representation with $*_2$ and $*_4$ both
zero. Write $W_{G,\p}$ for $V_{G,\p}/T_{G,\p}$.

We now show that $*_4$ gives us a non-zero class in
$\cohom_{f}^1(\rat, W_{f,\p}(1-k))$. Note that the fact that
$\overline{\rho}_{G,\p}$ is a homomorphism gives that $*_4$
necessarily gives a cohomology class in $\cohom^1(\rat,
W_{f,\p}(1-k)[\varpi])$.

First we suppose we are in the situation where $*_3 = 0$. Our
first step is to show that the quotient extension
\begin{equation}\label{eqn:quotientext1}
\mat{\overline{\rho}_{f,\p}}{*_4}{0}{\omega^{k-1}}
\end{equation}
is not split.  Suppose it is split.  Then by Proposition
\ref{thm:ribettypethm} we know that the extension
\begin{equation}
\mat{\omega^{k-2}}{*_2}{0}{\omega^{k-1}}
\end{equation}
cannot be split as well.  We show this gives a contradiction by
showing it gives a non-trivial quotient of the
$\omega^{-1}$-isotypical piece of the $p$-part of the class group
of $\rat(\mu_{p})$. However, Herbrand's Theorem (see for example,
\cite{washington}, Theorem 6.17) says that we must then have $p
\mid B_{2} = \frac{1}{30}$, which clearly cannot happen.

Consider the non-split representation
\begin{equation*}
\overline{\rho} = \mat{\omega^{-1}}{h}{0}{1}
\end{equation*}
which arises from twisting the non-split representation
\begin{equation*}
\mat{\omega^{k-2}}{*_2}{0}{\omega^{k-1}}
\end{equation*}
by $\omega^{1-k}$. Note that $\overline{\rho}$ is unramified away
from $p$ because $\overline{\rho}_{G,\p}$ is unramified away from
$p$.

We claim that this representation gives us a non-trivial finite
unramified abelian $p$-extension $K/\rat(\mu_{p})$ with the action
of $\Gal(K/\rat)$ on $\Gal(K/\rat(\mu_{p}))$ given by
$\omega^{-1}$.

Note that $\rat(\mu_{p}) = \overline{\rat}^{\ker \omega^{-1}}$, so
when we restrict $\overline{\rho}$ to
$\Gal(\overline{\rat}/\rat(\mu_{p}))$ we get
\begin{equation*}
\overline{\rho} \mid_{\Gal(\overline{\rat}/\rat(\mu_{p}))} =
\mat{1}{h}{0}{1},
\end{equation*}
i.e., we get a non-trivial homomorphism $h:
\Gal(\overline{\rat}/\rat(\mu_{p})) \rightarrow \mathbb{F}$ where $\mathbb{F}$ is a 
finite field of characteristic $p$. Set
$K = \rat(h) = \overline{\rat}^{\ker h}$, the splitting field of
$h$.

The fact that $\Gal(K/\rat(\mu_{p}))$ is abelian of $p$-power
order follows from the fact that
\begin{eqnarray*}\Gal(K/\rat(\mu_{p})) &\cong&
\Gal(\overline{\rat}/\rat(\mu_{p}))/\Gal(\overline{\rat}/\rat(h))
\\
&=& \Gal(\overline{\rat}/\rat(\mu_{p}))/ \ker h \\
&\cong& \image(h)
\end{eqnarray*}
and $\image(h)$ is a subgroup of $\mathbb{F}$, which is of
$p$-power order.  The fact that $K/\rat(\mu_{p})$ is unramified
away from $p$ also follows easily from the fact that
$\overline{\rho}$ is unramified away from $p$.  This shows that
$h(I_{\ell}) = 0$ for all $\ell \neq p$. In particular,
$h(I_{\ell}(K/\rat(\mu_{p}))) = 0$ for all $\ell \neq p$.  Since
we have the isomorphism above to a subgroup of $\mathbb{F}$, it
must be that $I_{\ell}(K/\rat(\mu_{p})) = 1$ for all $\ell \neq
p$.

The fact that $\Gal(K/\rat)$ acts on $\Gal(K/\rat(\mu_{p}))$ via
$\omega^{-1}$ follows from the fact that for $\sigma \in
\Gal(K/\rat(\mu_{p}))$ and $g \in \Gal(K/\rat)$, we have
\begin{equation*}
\overline{\rho}(g \sigma g^{-1}) \\
= \overline{\rho}(g) \overline{\rho}(\sigma)
\overline{\rho}(g^{-1}),
\end{equation*}
i.e., we have
\begin{equation*}
h(g \sigma g^{-1}) = \omega^{-1}(g) h(\sigma).
\end{equation*}

Our next step is to show that the extension $K/\rat(\mu_{p})$ that
we have constructed is actually unramified at $p$. We have that
$h|_{D_{p}} \in \cohom^1(\rat_{p}, \mathbb{F}(-1))$. Therefore, we
have that $h$ gives an extension $X$ of
$\ROI/\varpi \ROI$ by $\mathbb{F}(-1)$:\\

\begin{figure}[h!]
\centerline{\xymatrix{ 0 \ar[r]& \mathbb{F}(-1) \ar[r]& X \ar[r]&
\ROI/\varpi\ROI \ar[r]& 0.}}
\end{figure}

\noindent Applying Lemma \ref{thm:crystallineshort} and Lemma
\ref{thm:diamondflachguo2} we have that $h|_{D_{p}} \in
\cohom^1_{f}(\rat_{p}, \mathbb{F}(-1))$. A calculation in
\cite{blochkato} shows that $\cohom_{f}^1(\rat_{p}, E(-1)) = 0$
where $E$ is the field of definition for $\rho_{G,\p}$.  Actually,
it is shown that $\cohom_{f}^1(\rat_{p}, \rat_{p}(r))= 0$ for
every $r<0$; this implies
$\cohom_{f}^1(\rat_{p}, E(-1)) = 0$ since $E$ is a finite
extension (\cite{blochkato}, Example 3.9). Since we define
$\cohom_{f}^1(\rat_{p}, E/\ROI(-1))$ to be the image of the
$\cohom_{f}^1(\rat_{p}, E(-1))$, we have $\cohom_{f}^1(\rat_{p},
E/\ROI(-1)) = 0$.  Since $h|_{D_{p}} \in \cohom_{f}^1(\rat_{p},
\mathbb{F}(-1))$, Proposition \ref{thm:selmersummary} gives that
$h|_{D_{p}} \in \cohom^1_{f}(\rat_{p},E/\ROI(-1))$ and hence is 0.
Thus we have that $h$ vanishes on the entire decomposition group
$D_{p}$; in particular, it must be unramified at $p$ as
claimed.

Therefore, we have an unramified extension $K$ of $\rat(\mu_{p})$
that is of $p$-power order such that $\Gal(K/\rat)$ acts via
$\omega^{-1}$.  Let $C$ be the $p$-part of the class group of
$\rat(\mu_{p})$. Class field theory tells us that we have
\begin{equation*}
C/C^{p} \cong \Gal(F/\rat(\mu_{p}))
\end{equation*}
where $F$ is the maximal unramified elementary abelian
$p$-extension of $\rat(\mu_{p})$.  Therefore we have that
$\Gal(K/\rat(\mu_{p}))$ is a non-trivial subgroup of the
$\omega^{-1}$-isotypical piece of the $p$-part of the class group
of $\rat(\mu_{p})$, a contradiction as observed above.

Therefore, we must have that the quotient extension
\begin{equation*}
\mat{\overline{\rho}_{f,\p}}{*_4}{0}{\omega^{k-1}}
\end{equation*}
is not split if $*_3 =0$.

Now suppose that $*_1 = 0$.  Then the extension
\begin{equation*}
\mat{\omega^{k-2}}{*_2}{0}{\omega^{k-1}}
\end{equation*}
is a quotient extension and as above must necessarily be split.
Therefore again we get that the subextension
\begin{equation*}
\mat{\overline{\rho}_{f,\p}}{*_4}{0}{\omega^{k-1}}
\end{equation*}
cannot be split.

It remains to show that $*_4$ actually lies in $\cohom_{f}^1(\rat,
W_{f,\p}(1-k))$ since we have shown it is not zero. Write $h =
*_4$ to ease notation. As noted above, we have that $h$ gives a
non-zero class in $\cohom^1(\rat,W_{f,\p}(1-k)[\varpi])$.  Recall that in the previous
section we showed that $\cohom^1(\rat,W_{f,\p}(1-k)[\varpi])$ injects in $\cohom^1(\rat,
W_{f,\p}(1-k))$.  Therefore, we have that $h$ gives a non-zero class in $\cohom^1(\rat, W_{f,\p}(1-k))$.
 It remains to show that $h|_{D_{\ell}} \in
\cohom^1_{\ur}(\rat_{\ell}, W_{f,\p}(1-k))$ for each $\ell \neq p$
and $h|_{D_{p}} \in \cohom^1_{f}(\rat_{p}, W_{f,\p}(1-k))$. The
fact that $h|_{D_{\ell}} \in \cohom^1_{\ur}(\rat_{\ell},
W_{f,\p}(1-k)[\varpi])$ for $\ell \neq p$ is clear from the fact
that $\rho_{G,\p}$ is unramified away from $p$.  Therefore, we can
appeal to Proposition \ref{thm:diamondflachguo} to obtain that $h
\in \cohom^1_{\ur}(\rat_{\ell}, W_{f,\p}(1-k))$ for all $\ell \neq
p$.

The case at $p$ is easily handled by appealing to our work in the
previous section. Since $h|_{D_{p}} \in
\cohom^1(\rat_{p},W_{f,\p}(1-k)[\varpi])$, we get an extension $X$
of
$\ROI/\varpi \ROI$ by $W_{f,\p}(1-k)[\varpi]$:\\

\begin{figure}[h!]
\centerline{\xymatrix{ 0 \ar[r]& W_{f,\p}(1-k)[\varpi] \ar[r]& X
\ar[r]& \ROI/\varpi \ROI \ar[r]& 0. }}
\end{figure}

\noindent Appealing to Lemma \ref{thm:crystallineshort} and Lemma
\ref{thm:diamondflachguo2} we have that  $h|_{D_{p}}$ lies in
\linebreak$\cohom^1_{f}(\rat_{p}, W_{f,\p}(1-k)[\varpi])$ as
desired. Proposition \ref{thm:selmersummary} gives that
$h|_{D_{p}}$ lies in \linebreak $\cohom_{f}^1(\rat_{p},
W_{f,\p}(1-k))$.

Therefore, we have that $h$ is a non-zero torsion element of
$\cohom^1(\rat,W_{f,\p}(1-k))$ that lies in
$\cohom^1_{f}(\rat_{\ell}, W_{f,\p}(1-k))$ for every $\ell$.
Applying Proposition \ref{thm:selmersummary} to $h$ we have that
$h$ is a non-zero $\varpi$-torsion element of
$\cohom_{f}^1(\rat,W_{f,\p}(1-k))$. Therefore, it must be that $p
\mid \# \cohom_{f}^1(\rat,W_{f,\p}(1-k))$.  We summarize with the
following theorem.

\begin{theorem}\label{thm:finaltheorem} Let $k>3$ be an integer and $p>2k-2$ a
prime. Let \linebreak$f \in S_{2k-2}(\SL_2(\inte), \ROI)$ be a
newform with real Fourier coefficients so that
$\rho_{\mathfrak{m}_{f}}$ is irreducible and Conjecture
\ref{conj:heckeop} holds (for instance, if $f$ is ordinary at
$p$). Let
\begin{equation*}
\varpi^{m} \mid L_{\alg}(k,f)
\end{equation*}
with $m \geq 1$. If there exists an integer $N>1$, a fundamental discriminant $D$,
 and a Dirichlet character
$\chi$ of conductor $N$ so that
$(-1)^{k-1}D>0$, $\chi_{D}(-1)=-1$, $p \nmid ND[\Gamma_2:\Gamma_0^2(N)]$, and 
\begin{equation*}
\varpi^{n} \parallel L^{\Sigma}(3-k,\chi)L_{\alg}(k-1,f,\chi_{D})
L_{\alg}(1,f,\chi)L_{\alg}(2,f,\chi)
\end{equation*}
with $n < m$, then
\begin{equation*}
p \mid \# \cohom_{f}^1(\rat, W_{f,\p}(1-k)).
\end{equation*}
\end{theorem}

\section{Numerical Example}\label{sec:numericalexample}

In this concluding section we provide a numerical example of
Theorem \ref{thm:finaltheorem}.  We used the computer software
MAGMA, Stein's Modular Forms Database (\cite{stein}), and
Dokchitser's PARI program ComputeL (\cite{computel}).

Let $p = 516223$. We consider level 1 and weight 54 newforms in
$S_{54}(\SL_2(\inte))$.  There is one Galois conjugacy class of
such newforms, consisting of four newforms which we label $f_1,
f_2, f_3, f_4$.  Using the software Stein wrote for MAGMA we find
that
\begin{equation}\label{eqn:516223}
p \mid \prod_{i=1}^{4} L_{\alg}(28,f_{i}).
\end{equation}
The $q$-expansions of each $f_{i}$ are defined over a number field
$K_{i}$.  Appealing to MAGMA again we find each $K_{i}$ is
generated by a root of
\begin{eqnarray*}
\lefteqn{g(x) = x^4 + 68476320x^3 - 19584715019010048x^2}\\ & & \mbox{} -
10833127246634489297121280x\\ & & \mbox{} +
39446133467662904714689328971776.
\end{eqnarray*}
Let $\alpha_1, \alpha_2, \alpha_3, \alpha_4$ be the roots of
$g(x)$. Note that two of the $\alpha_{i}$ are real and the other
two are a complex conjugate pair. Relabelling the $f_{i}$ if
necessary, we may assume $K_{i} = \rat(\alpha_{i})$.  Let
$\ROI_{K_{i}}$ be the ring of integers of $K_{i}$.  Note that
$L_{\alg}(28,f_{i}) \in \ROI_{K_{i}}$ for each $i$. Therefore,
using Equation \ref{eqn:516223} we see that there exists $j \in
\{1,2,3,4\}$ and a prime $\wp_{j}\subset \ROI_{K_{j}}$ over $p$ so
that $\wp_{j} \mid L_{\alg}(28,f_{j})$. Since the $f_{i}$ are all
Galois conjugate, there is a conjugate prime $\wp_{i}\subset
\ROI_{K_{i}}$ over $p$ for each $i \in \{1,2,3,4\}$ so that
$\wp_{i} \mid L_{\alg}(28,f_{i})$.

Let $\chi = \chi_{-3}$ where we define $\chi_{-3}$  as in \cite{shimurahalfintegral}.  It is easy to
check that this $\chi$ and $D = -3$ satisfy the conditions of
Theorem \ref{thm:finaltheorem}. Using MAGMA we find that
\begin{equation*}
p \nmid \prod_{i=1}^{4} L_{\alg}(j,f_{i},\chi),
\end{equation*}
for $j=1,2$ and
\begin{equation*}
p \nmid \prod_{i=1}^{4} L_{\alg}(27,f_{i}, \chi_{D}).
\end{equation*}
We use ComputeL to show that
\begin{equation*}
p \nmid L^{(3)}(-25,\chi).
\end{equation*}
In particular, this shows we satisfy the divisibility hypotheses of Theorem \ref{thm:finaltheorem}
for $m=1$ and $n=0$.

Let $F_{i} = K_{i,\wp_{i}}$ with ring of integers $\ROI_{i}$ and
uniformizer $\varpi_{i}$. Set $\mathbb{F}_{i} =
\ROI_{i}/\varpi_{i} = \mathbb{F}_{p}[\overline{\alpha}_{i}]$ where
$\overline{\alpha}_{i} = \alpha_{i} (\modu \wp_{i})$. Let
$\rho_{i}: \absgal \rightarrow \GL_2(\ROI_{i})$ be the Galois
representation associated to $f_{i}$. This representation is
unramified away from $p$ and crystalline at $p$. Let
$\overline{\rho}_{i} : \absgal \rightarrow \GL_2(\mathbb{F}_{i})$
be the residual representation. Suppose that $\overline{\rho}_{i}$
is reducible.  Standard arguments show that $\overline{\rho}_{i}$
is non-split and we can write
\begin{equation*}
\overline{\rho}_{i} = \mat{\varphi}{*}{0}{\psi}
\end{equation*}
with $* \neq 0$ (see \cite{ribet1}).  Let $\omega: \absgal
\rightarrow \mathbb{F}_{p}^{\times}$ be the mod $p$ cyclotomic
character. Since $\varphi \psi = \omega^{53}$ and $\varphi$ and
$\psi$ are necessarily unramified away from $p$ and of order prime
to $p$, we can write $\varphi = \omega^{a}$ and $\psi =
\omega^{b}$ with $0 \leq a < b < p-1$, and $a + b = 53$ or $a+b = p
- 1 + 53$.  Arguing as in the previous section where we proved
that
\begin{equation*}
\mat{\overline{\rho}_{f,\p}}{*_4}{0}{\omega^{k-1}}
\end{equation*}
cannot be split, we have that $*$ gives a non-zero cocycle class
in $\cohom^1(\rat, \mathbb{F}_{i}(a-b))$ since $a-b <0$. As
before, this shows that we must have that $p$ divides the class
number of $\rat(\mu_{p})$, i.e., $p \mid B_{b-a+1}$ where we
recall that $B_{n}$ is the $n^{\text{th}}$ Bernoulli number
(\cite{washington}, Theorem 6.17). Appealing to the tables of
Buhler (\cite{buhler}), we see that the only Bernoulli number that
$516223$ divides is $B_{451304}$. Therefore, we must have $b-a+1 =
451304$, which in turn implies that $a+b = p - 1 + 53$ since
necessarily $a >0$. Solving this system of equations for $a$ and
$b$ we get $a = 32486$ and $b = 483789$. Observe that we have
\begin{eqnarray*}
\Tr(\overline{\rho}_{i}(\Frob_{2})) &=& 2^{a} + 2^{b} (\modu p)\\
        &=& 258573 (\modu p).
\end{eqnarray*}
Using Stein's tables we see that
$\Tr(\overline{\rho}_{i}(\Frob_2)) = \alpha_{i}$, so we must have
that $\overline{\alpha}_{i} \equiv 258573 (\modu \varpi)$. This
also shows that $\overline{\alpha}_{i}$ must belong to
$\mathbb{F}_{p}$ and so must be a root of one of the linear
factors of $g(x)$ modulo $p$. Using Maple to compute the linear
roots of $g(x)$ modulo $p$ we find that they are $287487$ and
$85284$, neither of which is congruent to $258573$ modulo $p$.
This provides a contradiction so we may conclude that
$\overline{\rho}_{i}$ is irreducible.

Due to the size of the prime under consideration, it was not
possible with the computer we used to compute the $p^{\text{th}}$
Fourier coefficients of the $f_{i}$ to check ordinarity. So,
instead we show that in this case the ordinarity assumption is not
necessary.  We do this by showing there are no congruences between
the $f_{i}$.  Let $E$ be a large number field containing all of
the $K_{i}$.  Let $i,j \in \{1,2,3,4\}$ with $i \neq j$. Let $\q$
be any prime of $E$ over $p$. As in Section \ref{sec:periods},
$f_{i}$ and $f_{j}$ each give a map from $\T_{\ROI_{E_{\q}}}$ to
$\ROI_{E_{\q}}$ given by $T(\ell) \mapsto a_{f_{i}}(\ell)$ and
$T(\ell) \mapsto a_{f_{j}}(\ell)$ respectively.  Let
$\mathfrak{m}_{i}$ and $\mathfrak{m}_{j}$ be the respective
maximal ideals defined as the inverse image of $\q$ under these
maps. (These are the maximal ideals associated to $f_{i}$ and
$f_{j}$ of $\T_{\ROI_{E_{\q}}}$ as in Section \ref{sec:periods}.)
There is a congruence between $f_{i}$ and $f_{j}$ modulo $\q$ if
and only if the maximal ideals $\mathfrak{m}_{i}$ and
$\mathfrak{m}_{j}$ are the same. This is equivalent to the
statement that
\begin{equation*}
a_{f_{i}}(\ell) \equiv a_{f_{j}}(\ell) (\modu \q)
\end{equation*}
for all $\ell \neq p$. In particular, looking at the case when
$\ell =2$, if a congruence exists between $f_{i}$ and $f_{j}$ we
have
\begin{equation*}
\q \mid (a_{f_{i}}(2) - a_{f_{j}}(2)),
\end{equation*}
i.e.,
\begin{equation*}
\q \mid (\alpha_{i} - \alpha_{j}).
\end{equation*}
Therefore we have that
\begin{equation*}
\Nm(\q) \mid \Nm(\alpha_{i} - \alpha_{j}).
\end{equation*}
The left hand side is a power of $p$ where as the right hand side
divides a power of the discriminant of $g(x)$, so that necessarily
$p$ divides the discriminant of $g(x)$.  Computing the
discriminant with Maple we find the prime factorization of the
discriminant,
\begin{eqnarray*}
\lefteqn{\disc(g(x)) = - 2^{48} 3^3 5^6 \cdot 11 \cdot 59 \cdot
15909926723 \cdot 4581597403}\\
&& \mbox{}\cdot
61912455248726091228769884731066259290896074682396020673553.
\end{eqnarray*}
Therefore we have that $p$ does not divide this discriminant.
Therefore we must have that there is no congruence modulo $\q$
between any of the $f_{j}$'s.  We can now appeal to the same
argument used in the proof of Lemma \ref{thm:congeigenform} to
conclude that there exists a Hecke operator $t$ so that $t \cdot
f_{i} = u \cdot \frac{\langle f_{i}, f_{i}
\rangle}{\Omega_{f_{i}}^{+}\Omega_{f_{i}}^{-}} f_{i}$ and $t \cdot
f_{j} = 0$ for $j \neq i$.  In this way we have avoided needing to
check the ordinarity of each $f_{j}$ to get the existence of the
Hecke operator conjectured in Conjecture \ref{conj:heckeop}.

If we choose $f_{i}$ to be one of the two newforms with real
Fourier coefficients, then we satisfy all of the hypotheses of
Theorem \ref{thm:finaltheorem} and so obtain the result that
\begin{equation*}
516223 \mid \# \cohom^1_{f}(\rat, W_{f_{i},\wp_{i}}(-27)).
\end{equation*}

\bibliographystyle{alpha}

\end{document}